\documentclass{amsart}

\usepackage{latexsym,amsthm,amssymb}

\usepackage[leqno]{amsmath}

\usepackage{mathrsfs}

\usepackage[all]{xy}

\usepackage{color}
\usepackage{colordvi}

\newtheorem{theorem}{Theorem}[section]
\newtheorem{lemma}[theorem]{Lemma}

\newtheorem{cor}[theorem]{Corollary}
\newtheorem{remark}[theorem]{Remark}

\newtheorem{observation}[theorem]{Observation}

\newtheorem{nota}[theorem]{Notation}

\numberwithin{equation}{theorem}

\newenvironment{Proofof}
{\noindent \emph{Proof of}}
{\hfill $\square$ \medskip}

\renewcommand{\mathcal}{\mathscr}

\renewcommand{\mathbb}{\mathbf}

\usepackage{verbatim}

\title{On the bicanonical morphism of quadruple Galois canonical covers}

\author{Francisco Javier Gallego
and Bangere P. Purnaprajna}

\address{Departamento de \'Algebra, Universidad Complutense de Madrid}
\email{gallego@mat.ucm.es}
\address{Department of Mathematics, University of Kansas}
\email{purna@math.ku.edu}
\thanks{The first author was partly supported by Spanish Goverment grant MTM2006-04785 and by Complutense grant PR27/05-13876 and is part of Complutense Research group 910772. He also thanks the Department of Mathematics of the University of Kansas for its hospitality. The second author thanks the General Research Fund of Kansas for partly supporting this research project. He also thanks the Department of Algebra of the Universidad Complutense de Madrid for its hospitality.}


\begin{document}
\begin{abstract}
In this article we study the bicanonical map $\varphi_2$ of quadruple Galois canonical covers $X$ of surfaces of minimal degree. We show that $\varphi_2$ has diverse behavior and exhibit most of the complexities that are possible for a bicanonical map of surfaces of general type, depending on the type of $X$.  There are cases in which $\varphi_2$ is an embedding, and if so happens, $\varphi_2$ embeds $X$ as a projectively normal variety, and cases in which $\varphi_2$ is not an embedding. If the latter, $\varphi_2$ is finite of degree $1$, $2$ or $4$. We also study the canonical ring of $X$, proving that it is generated in degree less than or equal to $3$ and finding the number of generators in each degree. For generators of degree $2$ we find a nice general formula which holds for canonical covers of arbitrary degrees. We show that this formula depends only on the geometric and the arithmetic genus of $X$. 
\end{abstract}
\maketitle
\section*{Introduction}

Canonical covers of surfaces of minimal degree have a ubiquitous presence in diverse contexts in 
the geometry of surfaces and threefolds. For example they appear in the classification of surfaces of 
general type 
with small $c_1^2$ as shown in the work of Horikawa,
and play an important role in mapping the geography of surfaces of
general type.  They appear as unavoidable boundary cases in the study of linear series on 
Calabi--Yau threefolds
as the works of Beltrametti and Szemberg (see ~\cite{BS}), Oguiso and Peternell (see~\cite{OP}) and the authors (see ~\cite{GPCalabi})
 show, and in the study of the canonical ring of a variety
 of general type as can be seen in the article of Green~\cite{Gr}.  
They are also a useful source in constructing new examples of
surfaces of general type. 
\medskip

Double and triple canonical covers were classified by Horikawa (see ~\cite{Ho}) and  Konno 
(see \cite{Kon}).
In~\cite{GPsing} and~\cite{GPsmooth}, the authors classified Galois canonical covers of degree $4$.
The classification showed that quadruple canonical covers behave quite differently from canonical covers
of all other degrees;  for instance, quadruple canonical covers are the only covers 
that admit 
families with unbounded geometric genus and families with unbounded irregularity. 
Hence, the geography of Chern numbers of quadruple canonical covers is much more complex.

\medskip

One of the much studied objects for surfaces of general type is the bicanonical map. In this article 
we prove results on the bicanonical morphisms of quadruple canonical covers of surfaces of minimal degree. 
Our results show that the behavior of these bicanonical maps 
is quite generic, that is, it exhibits all the diversities and complexities that are possible for a bicanonical map of a surface 
of general type. 
All of this is amply illustrated in the following theorem:  

\begin{theorem}
Let $\varphi: X \longrightarrow W$ be  a quadruple Galois canonical cover of a surface $W$ of minimal degree. 
Then the bicanonical map $\varphi_2$ of $X$ is
\begin{enumerate}
 \item[1)]  a morphism which embedds $X$ as a projectively normal variety if $X$ is of Type 1, 2, 3, 4, 5.2 or 6.2 (see Theorem~\ref{quad12} for the description of each type); 
\item[2)] a birational morphism but not an embedding if  $X$ is of Type 9, 10, 11 or 12; 
\item[3)] a finite morphism of degree $2$ if $X$ is of Type 5.1, 6.1, 7 or 8.2; 
\item[4)] a finite morphism of degree $4$ if $X$ is of Type 8.1. 
\end{enumerate}
\end{theorem}

The diverse behavior of the bicanonical maps exhibited in the above theorem is not seen in canonical covers of lower degrees and 
conjecturally does not happen for covers of all other degrees.

\medskip

In Section 2, we deal with those types of quadruple Galois canonical covers
 for which $2K_X$ is very ample and embeds $X$ as a projectively normal surface. Other than 
covers of the projective plane, the manifestation of such behavior can be seen in surfaces $X$ with unbounded $p_g$ but with bounded $q$. 
However, in the case of families with unbounded $q$, $2K_X$ is not very ample, even though the image is a projectively normal variety.
 One of the
frequently used techniques to show that a $2$--Veronese  subring of a graded ring is generated in degree less than or equal to two is to first show
 that the ambient graded ring is generated in degree less than or equal to two.  Formally one can construct graded rings in commutative algebra 
where this does not happen and yet the $2$--Veronese  subring is generated in degree $1$.  In this context, the quadruple covers of types 1, 2, 3, 4, 5.2 and 6.2  provide natural
examples where the canonical ring is not generated in degree less than or equal to two, yet the $2$--Veronese  subrings are generated in degree less 
or equal one,
thereby showing the normal generation of the bicanonical map.

\medskip

Section 3 deals with those quadruple covers for which bicanonical maps are finite. The types for which the bicanonical maps are not birational 
have unbounded $p_g$ and $q$. But the results of Xiao (see~\cite[Theorems 1, 2 and 3]{Xiao}) say that if the bicanonical maps are not birational then 
they form bounded families with respect to $p_g$ 
unless the surfaces possess a genus two pencil.
In Section 3, 
we explicitly exhibit the genus two pencils that are indeed fibrations.
So in addition to being a genus $3$ fibration over 
$\mathbf P^1$, the families of types 5.1 and 6.1 are also genus two fibrations over an elliptic curve, and the families of type 
7 are also genus 2 fibrations over a curve of genus $m$, where $m$ takes on unbounded values.

\medskip

Section 4 deals with the Types 9, 10, 11 and 12 where the image of the canonical map $W$ is a singular surface of minimal degree. It was shown in 
\cite{GPsing} that in such a case $p_g\leq 4$ and $q=0$. The behavior of the bicanonical map is very interesting for these types of surfaces: it is always 
birational but never an embedding. And, more interestingly, it is birational is two different ways: for the types 9, 10 or 12, $|2K_X|$ does not 
separate directions at the unique point $x \in \varphi^{-1}\{w\}$ (and is an isomorphism outside $x$). For type 11, $|2K_X|$ does not separate the two points $x_1$ and $x_2$ of 
$\varphi^{-1}\{w\}$, although ${\varphi_2}$ is locally an embedding at both of them (and it is an isomorphism outside $x_1$ and $x_2$).  We prove these results by proving a non--vanishing 
theorem for certain ideal sheaves.  To accomplish this, we construct explicit factorizations of some birational maps which are in half the cases 
crepant and in the other half not and handle the question of the non--vanishing result on a suitably constructed desingularization of $X$. In all of this, the algebra 
structure of the map $p$ in~\eqref{CD} (the so--called desingularized diagram), which was precisely described in \cite{GPsing}, also plays an important role. 

\medskip

The canonical ring of a surface of general type and the degrees of its generators have attracted  interest among geometers for various reasons. One such reason 
is its applications to the study of Calabi--Yau threefolds.
For example, results on ring generation are used in determining the very ampleness  of line bundles on
Calabi--Yau threefolds in the article~\cite{GPCalabi} and have provided motivation in the construction of new examples of Calabi--Yau threefolds as can be seen in the work of Casnati~\cite{Cas}.
In Section 5 we prove a general result, Theorem~\ref{gen.in.deg2},
that gives a nice formula for the number of generators in degree $2$ of 
the canonical ring of canonical covers, of arbitrary degree, of surfaces of minimal degree. The formula shows that the number of generators 
in degree $2$ only depends on the geometric and arithmetic genus of $X$. 
Theorem~\ref{canonical.ring}, which determines the generators of the canonical ring of quadruple covers, shows that there is no such formula
for generators of degree $3$ of the canonical ring of $X$, if $X$ is an irregular surface of general type. In fact Theorem~\ref{canonical.ring} 
shows that this number depends on the algebra structure of $\varphi$. 

\medskip

\noindent {\emph{Acknowledgements}:} We thank Dale Cutkosky and N. Mohan Kumar 
for some very valuable discussions and comments.

\section{Preliminaries and notation}

We introduce the notation and the basic definitions that we will use throughout the article: 

\medskip

\noindent{\bf Convention.}
 We work over an algebraically closed field of
 characteristic $0$.

\begin{nota}\label{notations}
{\rm {Throughout this article, unless otherwise stated, we will make the following assumptions and use the following notation:

\medskip
\begin{enumerate}
\item[(1)] 
$W$ will be an embedded projective
   algebraic surface of
 minimal degree, i.e., a surface such that deg$W=$codim$W+1$.
\medskip
\item[(2)]
$X$ will be a projective
 algebraic normal
 surface
with at worst canonical singularities (that is, $X$ is smooth or has rational double points; see ~\cite{Ba} for details about rational double points), whose canonical divisor
$K_{X}$ is ample and base--point--free. 
\medskip
\item[(3)] We will denote the canonical map of $X$ as  $\varphi$. Note that, by (2), $\varphi$ is in fact a finite morphism.
\medskip
\item[(4)] We will assume  $\varphi$ to be a  Galois morphism of degree $4$ whose image is a surface $W$ of minimal degree, that is, $\varphi: X \longrightarrow W$ will be a {\it quadruple Galois canonical cover of a surface $W$ of minimal degree}. 
\medskip
\item[(5)] 
We will denote the bicanonical map of $X$ as $\varphi_2$.

\end{enumerate}
\medskip

\noindent  We also recall the following standard notation, that will also be used throughout the article:

\medskip

\begin{enumerate}

\item[(6)] By $\mathbf F_e$ we denote the Hirzebruch surface whose
 minimal section have self-intersection $-e$. If $e >0$ let $C_0$
 denote the minimal section of $\mathbf F_e$ and let $f$ be one of the
 fibers of $\mathbf F_e$. If $e=0$, $C_0$ will be a fiber of one of the
 families of lines and $f$ will be a fiber of the other family of
 lines of $\mathbf F_0$.

\medskip

\item[(7)] If $a, b$ are integers such that $0 < a \leq b$, consider
 two disjoint linear subspaces $\mathbf P^a$ and $\mathbf P^b$ of $\mathbf
 P^{a+b+1}$. We denote
   by $S(a,b)$ the smooth rational normal scroll obtained by joining
   corresponding points of a rational normal curve in $\mathbf P^a$
   and a rational normal curve of $\mathbf P^b$.
Recall that $S(a,b)$ is the
image of $\mathbf F_e$ by the embedding induced by the complete linear series
$|C_0+mf|$, with $a=m-e$, $b=m$ and $m > e+1$.

\medskip

\item[] If $a=b$, the linear series
$|mC_0+f|$ also gives a minimal degree embedding of $\mathbf F_0$,
equivalent to the previous one by
the automorphism of $\mathbf P^1 \times \mathbf P^1=\mathbf F_0$ swapping the
factors. In this case {\it {our convention will always be to choose $C_0$ and
$f$ so that, when $W=S(m,m)$, $W$ is
embedded by $|C_0+mf|$.}}

\medskip

\item[] If in addition $m=1$, $C_0$ and $f$ are indistinguisable in both $\mathbf F_0$
and $S(1,1)$, so, {\it {in such a case, for us $C_0$ will denote the fiber of any of the
families of lines of $\mathbf F_0$ and $f$ will denote the fiber of the other family}}.

\medskip

\item[(8)] If $b$ is an integer, $b > 1$, consider
 a linear subspace $\mathbf P^b$ of $\mathbf
 P^{b+1}$. We denote
   by $S(0,b)$ the cone in $\mathbf
 P^{b+1}$ over a rational normal curve of $\mathbf P^b$.
Recall that $S(0,b)$ is the
image of $\mathbf F_e$ by the morphism induced by the complete linear series
$|C_0+mf|$, with $b=m=e$ and hence $e > 1$.
\end{enumerate}}}
\end{nota}

\begin{remark}\label{Galrem}
\rm{If $X \overset{\varphi} \longrightarrow W$ is a Galois cover and $W$ is smooth, then $\varphi$ is flat.}
\end{remark}

We recall now the main features of the classification of $\varphi$, which was obtained in~\cite[Theorem 0.1]{GPsmooth} and in the main theorem of~\cite{GPsing}. According to this classification $\varphi$ falls into several different types described in the tables of the theorem below. We will refer to these types throughout the article. 

\medskip

\begin{theorem}\label{quad12}
Let $\varphi: X \longrightarrow W$ be as in Notation~\ref{notations} (i.e., let $\varphi$ be a quadruple Galois canonical cover of a surface $W$ of minimal degree). 

\begin{enumerate} 
\item[\bf 1)]  If $W$ is {\rm smooth}, then $W$ is either linear $\mathbf P^2$  or a smooth Hirzebruch surface $\mathbf F_e$,
with $0 \leq e \leq 2$, embedded by $|H|$, where $H=C_0+mf$ ($m \geq e+1$). 
Recall that the Galois group $G$ of $\varphi$ is either $\mathbf Z_4$ or $\mathbf Z_2^{\oplus 2}$.

\smallskip

\noindent If $G=\mathbf Z_4$, then
$\varphi$ is the composition of two double covers
$X_1 \overset{p_1} \longrightarrow W$ branched along a divisor $D_2$
and $X \overset{p_2} \longrightarrow X_1$, branched along the
ramification of $p_1$ and $p_1^*D_1$, where $D_1$ is a divisor on $W$.

\smallskip

\noindent If $G= \mathbf Z_2^{\oplus 2}$, then
$X$ is the fiber product over $W$ of two
double covers of $W$ branched along divisors $D_1$ and $D_2$ and $\varphi$
is the natural morphism from the fiber product to $W$.

\smallskip

\noindent Moreover $\varphi$ fits into one of the following types, determined by these characteristics:

\bigskip

\centerline{\vbox{\tabskip=0pt \offinterlineskip
\def\tablerule{\noalign{\hrule}}
\halign to 16truecm
{\strut
#& \vrule#
\tabskip=0em plus 3em
&
\hskip .3truecm \hfil# \hskip .1truecm
& \vrule #
& \hskip .2truecm \hfil #
\hfil  \hskip .1truecm & \vrule#&
\hskip .2truecm \hfil#\hfil \hskip .1truecm & \vrule#&
\hskip .2truecm \hfil#\hfil \hskip .1truecm & \vrule#&
\hskip .2truecm \hfil#\hfil \hskip .1truecm & \vrule#&
\hskip .2truecm \hfil#\hfil \hskip .1truecm &
\vrule#&
\hskip .25truecm \hfil# %\hskip .05truecm
\hfil & \vrule#\tabskip=0pt
\cr\tablerule
&&Type
&&%\omit\hidewidth
$W$ %\hidewidth
&& %\omit\hidewidth
$p_g(X)$ %\hidewidth
&&%\omit\hidewidth
$G$ %\hidewidth
&&%\omit\hidewidth
$D_1 \sim$ %\hidewidth
&&%\omit\hidewidth
 $D_2 \sim$  %\hidewidth
&&%\omit\hidewidth
  $q(X)$  %hidewidth
%&&%\omit\hidewidth
%$c_1^2/c_2$ %\hidewidth
&\cr\tablerule
\cr
\tablerule
&&1
&& $\mathbf P^2$
&&$3$
&&  $\mathbf Z_4$
&&  conic
&&  quartic
&& $0$
%&&  $\frac{1}{11}$
&\cr
\tablerule
&&2
&&$\mathbf P^2$
&&$3$
&&$\mathbf Z_2^{\oplus 2}$
&&quartic
&&quartic
&&$0$
%&&$\frac{1}{11}$
&\cr
\tablerule
&&3
&&$S(m,m-e)$  %\ *$ \hskip -.2truecm
&&$2m-e+2$
&&$\mathbf Z_4$
&&$(2m-e+1)f$
&&$4C_0+(2e+2)f$
&&$0$
%&&$\frac{2m-e}{4m-2e+9}$
&\cr
\tablerule
&&4
&&$S(m,m-e)$ % \ *$ \hskip -.2truecm
&&$2m-e+2$
&&$\mathbf Z_2^{\oplus 2}$
&&$2C_0+(2m+2)f$
&&$4C_0+(2e+2)f$
&&$0$
%&&$\frac{2m-e}{4m-2e+9}$
&\cr
\tablerule
&&5.1
&&$S(1,1)$ %\hskip -.15truecm
&&$4$
&&$\mathbf Z_4$
&&$6f$
&&$4C_0$
&&$1$
%&&$\frac{m}{2m+3}$
&\cr
\tablerule
&&5.2
&&$S(m,m)$, ($m>1$) %\hskip -.15truecm
&&$2m+2$
&&$\mathbf Z_4$
&&$(2m+4)f$
&&$4C_0$
&&$1$
%&&$\frac{m}{2m+3}$
&\cr
\tablerule
&&6.1
&&$S(1,1)$ %\hskip -.15truecm
&&$4$
&&$\mathbf Z_2^{\oplus 2}$
&&$4C_0$
&&$2C_0+6f$
&&$1$  
&\cr
\tablerule
&&6.2
&&$S(m,m)$, ($m>1$)%\hskip -.15truecm
&&$2m+2$
&&$\mathbf Z_2^{\oplus 2}$
&&$2C_0+(2m+4)f$
&&$4C_0$
&&$1$
%%&&$\frac{m}{2m+3}$
&\cr
\tablerule
&&7
&&$S(m,m)$ %\hskip -.15truecm
&&$2m+2$
&&$\mathbf Z_2^{\oplus 2}$
&&$(2m+2)f$
&&$6C_0+2f$
&&$m$
%&&$\frac{2m}{m+9}$
&\cr
\tablerule
&&8.1
&&$S(1,1)$ %\hskip -.15truecm
&&$4$
&&$\mathbf Z_2^{\oplus 2}$
&&$6f$
&&$6C_0$
&&$4$
%&&$2$
&\cr
\tablerule
&&8.2
&&$S(m,m)$, ($m>1$)%\hskip -.15truecm
&&$2m+2$
&&$\mathbf Z_2^{\oplus 2}$
&&$(2m+4)f$
&&$6C_0$
&&$m+3$
%&&$2$
&\cr
\tablerule
\noalign{\smallskip}}}}

\medskip

\item[\bf 2)] If $W$ is {\rm not smooth}, then $W=S(0,2)$ (a quadric cone in $\mathbf P^3$), $X$ is regular  and
$\varphi: X \longrightarrow W$ fits into a commutative diagram (see~\cite{GPsing}):

\begin{equation}\label{CD}
\xymatrix@1{\overline X \ar[d]^p \ar[r]^{\overline \psi} &   X \ar[d]^{\varphi} \\
Y \ar[r]^\psi & W}, 
\end{equation}

where $\psi$ is the minimal desingularization of $W$, $\overline X$ is the normalization of
 the reduced part 
of $X \times_W Y$ and has at worst canonical singularities, and $p$ and $\overline \psi$ are
 induced by the projections from the
 fiber product onto each factor, $\overline \psi$ being the morphism from $\overline X$ to its canonical model $X$. The structure of $p$ is similar to the structure of $\varphi$ described above for $W$ smooth (see the main theorem of~\cite{GPsing} for details) and $\varphi$ can be classified in another four types, determined by the following properties: 
  
\bigskip

\centerline{\vbox{\tabskip=0pt \offinterlineskip
%\tabskip= .25 truecm
\def\tablerule{\noalign{\hrule}}
\halign to 16truecm
%{\valign to125pt}
{\strut
#& \vrule#
\tabskip=0em plus 3em
%\tabskip= .25 truecm
&
\hskip .3truecm \hfil# \hskip .1truecm
& \vrule #
& \hskip .2truecm \hfil #
\hfil  \hskip .1truecm & \vrule#&
\hskip .2truecm \hfil#\hfil \hskip .1truecm & \vrule#&
\hskip .2truecm \hfil#\hfil \hskip .1truecm & \vrule#&
%\hskip .2truecm \hfil#\hfil \hskip .1truecm & \vrule#&
\hskip .2truecm \hfil#\hfil \hskip .1truecm &
\vrule#&
\hskip .25truecm \hfil# %\hskip .05truecm
\hfil & \vrule#\tabskip=0pt
\cr\tablerule
&&Type
&&%\omit\hidewidth
$W$ %\hidewidth
&& %\omit\hidewidth
$p_g(X)$ %\hidewidth
&&%\omit\hidewidth
$G$ %\hidewidth
&&%\omit\hidewidth
$\overline \psi$ is%\hidewidth
%&&%\omit\hidewidth
 %$D_2 \sim$  %\hidewidth
&&%\omit\hidewidth
  $q(X)$  %hidewidth
%&&%\omit\hidewidth
%$c_1^2/c_2$ %\hidewidth
&\cr\tablerule
\cr
\tablerule
&&9
&& $S(0,2)$
&&$4$
&&  $\mathbf Z_2^{\oplus 2}$
&&  crepant
%&&  
&& $0$
%&&  
&\cr
\tablerule
&&10
&&$S(0,2)$
&&$4$
&&$\mathbf Z_4$
&&crepant
%&&
&&$0$
%&&
&\cr
\tablerule
&&11
&&$S(0,2)$  %\ *$ \hskip -.2truecm
&&$4$
&&$\mathbf Z_2^{\oplus 2}$
&&non crepant
%&&$
&&$0$
%&&$\frac{2m-e}{4m-2e+9}$
&\cr
\tablerule
&&12
&&$S(0,2)$ % \ *$ \hskip -.2truecm
&&$4$
&&$\mathbf Z_4$
&&non crepant
%&&
&&$0$
%&&
&\cr
\tablerule
\noalign{\smallskip}}}}
\end{enumerate}
\end{theorem}

\section{Normal generation of  the bicanonical bundle}\label{Normal.generation.2K_X}

In this section we study the types of quadruple Galois canonical covers $X$ for which $2K_X$ is normally generated, i.e., we find out for what $X$ the bicanonical morphism $\varphi_2$ embeds $X$ as a projectively normal variety. As we will see in Section~\ref{sing.section}, $\varphi_2$ is never an embedding if $W$ is singular, so throughout this section we assume $W$ \emph{to be smooth}. 

\medskip

\noindent Recall that by Remark~\ref{Galrem}, since we are assuming $W$ to be smooth, $\varphi$ is flat. Furthermore, the push down of $\mathcal O_X$ to $W$ splits as an $\mathcal O_W$--module as follows: 
\begin{equation}\label{split}
 \varphi_* \mathcal{O}_{X}=\mathcal{O}_{W}\oplus L _{1}^* \oplus
L _{2}^* \oplus L _{3}^* , 
\end{equation}
with $L_1, L_2$ and $L_3$ and $D_1$ and $D_2$ of Theorem~\ref{quad12} satisfying the following properties (see~\cite[Remark 3.1]{GPsmooth}, where the multiplicative structure that turns the second term of~\eqref{split} into an $\mathcal O_W$--algebra is also described): 
\begin{equation}\label{split2}
\begin{matrix}
&L_1 \otimes L_2 = L_3 &\\ 
\textrm{ if } G=\mathbf Z_2 \times \mathbf Z_2, & \textrm{ then } L_1^*=\mathcal O_W(-1/2D_2) \textrm{ and } & L_2^*=\mathcal O_W(-1/2D_1); \textrm{ and } \\ 
 \textrm{ if } G=\mathbf Z_4, \textrm{ then } & L_1^*=\mathcal O_W(-1/2D_1-1/4D_2) &  \textrm{ and } L_2^*=\mathcal O_W(-1/2D_2). 
\end{matrix}
\end{equation}
We will keep the notation introduced in~\eqref{split} and~\eqref{split2} for the remaining of the article. 

\medskip

\noindent To prove or disprove the normal generation of $2K_X$ we will look at multiplication maps of global sections of line bundles on $X$. To study these maps we will use the $\mathcal O_W$--algebra structure of $\varphi_*\mathcal O_X$ as the following lemma explains:

\begin{lemma}\label{surject.algstruct}
Let $A_{1},$ $A_{2}$ be
two line bundles on $W$ and let $M_1=\varphi^* A_{1}$ and
$M_2=\varphi^* A_{2}$ be their pull backs on $X$. 

Let%
\begin{equation*}
H^{0}(M_1)\otimes H^{0}(M_2)\overset{\beta}{\longrightarrow}H^{0}(M_1\otimes M_2)
\end{equation*}
be the multiplication map of global sections of $M$ and $N$ and let  
\begin{equation*}
\begin{matrix}
H^{0}(A_{1})\otimes H^{0}(A_{2}) & \overset{\beta_{1}}{\longrightarrow } & 
H^{0}(A_{1}\otimes A_{2}), \\
H^{0}(A_{1})\otimes H^{0}(A_{2}\otimes
L _{1}^* ) & \overset{\beta_{2}}{\longrightarrow } & H^{0}(A_{1}\otimes
A_{2}\otimes L _{1}^* ), \\  
H^{0}(A_{1})\otimes H^{0}(A_{2}
\otimes L _{2}^* ) & \overset{\beta_{3}}{\longrightarrow } & H^{0}  
(A_{1}\otimes A_{2}\otimes L _{2}^* ), \\
H^{0}(A_{1}%
\otimes L _{1}^* )\otimes H^{0}(A_{2}\otimes L _{2}^{\ast
}) & \overset{\beta_{4}}{\longrightarrow } & H^{0}(A_{1}\otimes A_{2}\otimes
L _{3}^* )
\end{matrix}
\end{equation*}
be multiplication maps of global sections of line bundles on $W$.

\begin{enumerate}
\item[1)] If $\beta_1$, $\beta_2$, $\beta_3$ and $\beta_4$ are surjective, so is $\beta$. 

\item[2)] 
If $H^{0}(A_{1}\otimes
A_{2}\otimes L _{3}^* )\neq 0$ but   $H^{0}%
(A_{i}\otimes L _{j}^* )=0$ for either 
\begin{enumerate}
\item[] $i=1$ and $j=1,2$,  or
\item[] $i=2$ and $j=1,2$, or
\item[] $i=1,2$ and $j=1$, or
\item[] $i=1,2$ and $j=2$,  
\end{enumerate}
then  $\beta$ is not surjective. 
\end{enumerate}
\end{lemma}

\begin{proof} 
We have by
projection formula
\begin{equation}\label{split.pushdown}
\begin{matrix}
H^{0}(\varphi_* M_i) =  H^{0}(A_{i})\oplus H^{0}(A_{i}\otimes L _{1}^{\ast
})\oplus H^{0}(A_{1}\otimes L _{2}^* )\oplus H^{0}(A_{1}%
\otimes L _{3}^* )  \textrm{ and } \\  
H^{0}(\varphi_*(M_1\otimes M_2))=  H^{0}(A_{1}\otimes A_{2})\oplus H^{0}%
(A_{1}\otimes A_{2}\otimes L _{1}^* )\oplus H^{0}(A_{1}\otimes
A_{2}\otimes L _{2}^* )\oplus H^{0}(A_{1}\otimes A_{2}%
\otimes L _{3}^* ). 
\end{matrix}
\end{equation}

The surjectivity
of $\beta$ is equivalent to the surjectivity of
\[
H^{0}(\varphi_* M)\otimes H^{0}(\varphi_* N)\overset{\beta^{\prime}}{\longrightarrow 
}H^{0}(\varphi_* (M\otimes N)). 
\]
The $\mathcal O_W$--algebra structure of $\varphi_*\mathcal O_X$ is given by a multiplication map
\begin{equation*}
 \varphi_*\mathcal O_X \otimes \varphi_*\mathcal O_X \longrightarrow \varphi_*\mathcal O_X
\end{equation*}
which splits in several summands according to~\eqref{split}, as explained in~\cite[Remark 3.1]{GPsmooth}. From them, we are interested in  
the following four:%

\begin{equation}\label{mult.struct}
\begin{matrix}
\mathcal{O}_{W}\otimes\mathcal{O}_{W} &  \overset{\simeq}\longrightarrow \mathcal{O}_{W}\\
\mathcal{O}_{W}\otimes L _{1}^*  &  \overset{\simeq}\longrightarrow
L _{1}^*\\
\mathcal{O}_{W}\otimes L _{2}^*  &  \overset{\simeq}\longrightarrow
L _{2}^*\\
L _{1}^*\otimes L _{2}^*  &  \overset{\simeq}\longrightarrow
L _{3}^*%
\end{matrix}
\end{equation}

The map $\beta'$ also splits according to~\eqref{split.pushdown}, so~\eqref{mult.struct} implies that if $\beta_1,\beta_2,\beta_3$ and $\beta_4$ surject, all the summands of $H^{0}(\varphi_*(M_1\otimes M_2))$ described in~\eqref{split.pushdown} are in the image of $\beta'$, so $\beta'$ and hence $\beta$ are surjective. On the other hand, 
\cite[Remark 3.1]{GPsmooth} also tells that the last map  of~\eqref{mult.struct}  is the only summand mapping to $L_3$.  
Thus, under the hypotheses of 2), $H^{0}(A_{1}\otimes A_{2}\otimes
L _{3}^* ) \neq 0$ but it is not in the image of $\beta'$, so $\beta'$ and hence $\beta$ are not surjective in this case. 
\end{proof}

To apply Lemma~\ref{surject.algstruct} in proving the normal generation of $2K_X$ we will need this easy but useful observation that helps to handle multiplication maps of global sections:

\begin{observation}\label{horace} Let $E$ and $L_{1},...,L_{n}$ be
coherent sheaves on a variety $X$. Consider the map $H^{0}(E)\otimes H^{0}%
(L_{1}\otimes\cdot\cdot\cdot\otimes L_{r})\overset{\psi}{\longrightarrow }%
H^{0}(E\otimes L_{1}\otimes\cdot\cdot\cdot\otimes L_{r})$ and the maps
\begin{align*}
&  H^{0}(E)\otimes H^{0}(L_{1})\overset{\alpha_{1}}{\longrightarrow }H^{0}(E\otimes
L_{1})\\
&  H^{0}(E\otimes L_{1})\otimes H^{0}(L_{2})\overset{\alpha_{2}}{\longrightarrow 
}H^{0}(E\otimes L_{1}\otimes L_{2})\\
&  \text{ \ \ \ \ \ \ \ \ \ \ \ \ \ \ \ \ \ \ \ \ \ \ \ \ \ \ \ }........\\
&  H^{0}(E\otimes L_{1}\otimes\cdot\cdot\cdot\otimes L_{r-1})\otimes
H^{0}(L_{r})\overset{\alpha_{r}}{\longrightarrow }H^{0}(E\otimes L_{1}\otimes
\cdot\cdot\cdot\otimes L_{r})
\end{align*}
If $\alpha_{1},...\alpha_{r}$ are surjective then $\psi$ is
surjective.
\end{observation}

Now we are ready to prove the normal generation of $2K_X$ for surfaces $X$ of Types 1, 2, 3, 4, 5.2 and 6.2. As a warm--up for the rest, we start with the simplest case, that is, when $W=\mathbf P^2$: 

\begin{theorem}\label{types12}
 
Let $\varphi:X\longrightarrow  W$ be a quadruple Galois canonical cover of  Type 1 or 2. Then $2K_{X}$ is
very ample and $|2K_{X}|$ embeds $X$ as a projectively normal variety.
\end{theorem}

\begin{proof} 
Recall that $K_X=\varphi^*\mathcal O_{\mathbf P^2}(1)$. Since $K_X$ is ample, the normal generation of $2K_X$ is equivalent to the surjectivity of 
\begin{equation}\label{projnorm}
 H^0(2K_X) \otimes H^0(2nK_X) \longrightarrow H^0((2n+2)K_X), 
\end{equation}
for all $n \geq 1$. We prove it in two steps. 

\smallskip

\noindent\emph{Step 1.} The 
first step is to show that
\begin{equation*}
H^{0}(2K_{X})\otimes H^{0}(2K_{X})\overset{\beta}{\longrightarrow }H^{0}%
(4K_{X}) %
\end{equation*}
surjects. We
know (see the structure of $\varphi$ as described in Theorem~\ref{quad12}) that, for both Types 1 and 2, 
\begin{equation*}
\varphi_* \mathcal{O}_{X}=\mathcal{O}_{\mathbf P^2}\oplus\mathcal{O}_{\mathbf P^2}(-2) \oplus \mathcal{O}_{\mathbf P^2}(-2)\oplus \mathcal{O}_{\mathbf P^2}
(-4), 
\end{equation*}
hence we have 
\begin{eqnarray*} 
H^{0}(\varphi_*2K_{X})&=&H^{0}(\mathcal{O}_{\mathbf P^2}(2))\oplus H^{0}%
(\mathcal{O}_{\mathbf P^2})\oplus H^{0}(\mathcal{O}_{\mathbf P^2}) \textrm{ \ \ and }\\
H^{0}(\varphi_*4K_{X})&=&H^{0}(\mathcal{O}_{\mathbf P^2}(4))\oplus H^{0}%
(\mathcal{O}_{\mathbf P^2}(2))\oplus H^{0}(\mathcal{O}_{\mathbf P^2}(2))  \oplus H^{0}(\mathcal{O}_{\mathbf P^2}). 
\end{eqnarray*} 
Then Lemma~\ref{surject.algstruct}, 1) tells that in order to prove the surjectivity of $\beta$ it is enough to show the surjectivity of the
following multiplication maps on $W$:%
\begin{eqnarray*}
 H^{0}(\mathcal{O}_{\mathbf P^2}(2))\otimes H^{0}(\mathcal{O}_{\mathbf P^2}(2))&\overset{\beta_{1}}{\longrightarrow 
}&H^{0}(\mathcal{O}_{\mathbf P^2}(4))\\
H^{0}(\mathcal{O}_{\mathbf P^2})\otimes H^{0}(\mathcal{O}_{\mathbf P^2}(2))&\overset{\beta_{2}}{\longrightarrow 
}&H^{0}(\mathcal{O}_{\mathbf P^2}(2))\\
H^{0}(\mathcal{O}_{\mathbf P^2})\otimes H^{0}(\mathcal{O}_{\mathbf P^2})&\overset{\beta_{4}}{\longrightarrow 
}&H^{0}(\mathcal{O}_{\mathbf P^2}).
\end{eqnarray*}
The surjectivity of $\beta_1$ follows from the projective normality of the Veronese surface and the surjectivity of $\beta_2$ and $\beta_4$ is trivial. 

\smallskip

\noindent\emph{Step 2.} To complete the proof, we need to show that the
multiplication map~\eqref{projnorm}
surjects for all
$n \geq 2$. In view of Observation~\ref{horace}, it is enough to show that
\begin{equation*} 
H^{0}(n'K_{X})\otimes H^{0}(K_{X})\longrightarrow 
H^{0}((n'+1)(K_{X})), 
\end{equation*}
for all $n' \geq 4$.  Since $K_X$ is ample and base--point--free, this 
follows from~\cite[p. 41, Theorem 2]{Mum} and the Kawamata--Viehweg vanishing theorem.
\end{proof}

Now we go on to study the normal generation of  $2K_X$ for those remaining cases  for which $W$ is smooth.  In these cases $W$ is a Hirzebruch surface, so in order to apply Lemma~\ref{surject.algstruct} we need to know first the surjectivity of certain multiplication maps of global sections of line bundles on Hirzebruch surfaces. This is done in the next lemma, where we give sufficient conditions for the surjectivity of such maps. 

\begin{lemma}\label{surjonHirz}
Let $W$ $=\mathbf{F}_{e}, e=0,1$ or $2.$ 
Let
$L_{1}=a_{1}C_{0}+b_{1}f$ and $L_{2}=a_{2}C_{0}+b_{2}f$ be two line bundles on $W$, with $a_{i}\geq0,$
$b_{i}\geq a_{i}e$.  Then the multiplication map
\begin{equation*}
H^{0}(L_{1})\otimes H^{0}(L_{2})\overset{\mu}{\longrightarrow }H^{0}(L_{1}%
+L_{2})
\end{equation*}
is surjective if in addition $L_1$ and $L_2$ satisfy one of the following conditions:
\begin{enumerate}
\item[(a)] $a_{1}\geq1, a_{2}=na_{1},
b_{2}=nb_{1}$ with $n\geq1$,  and $b_{1}>a_{1}e;$
\item[(b)] $a_{1}\geq1$  and  $a_{2}\geq2a_{1}-2+e$ if
 $e\geq1$  or  $a_{2}\geq2a_{1}-1$  if  $e=0$, 
  and
 $b_{2}-b_{1}\geq (a_2-a_1)e-1;$
\item[(c)] $a_{1}>0,a_{2}=0$;
\item[(d)] $a_{1}=a_{2}=1,b_{2}\geq b_{1}-1;$
\item[(e)] $W=\mathbf{F}_{0}.$ 
\end{enumerate}
\end{lemma}

\begin{proof}
Observe that $L_1$ and $L_2$ are
base--point--free because $a_{i}\geq0,$ $b_{i}\geq a_{i}e$. 

\smallskip

We start proving the lemma under the assumption that (a) is satisfied. Observe that $L_1$ is ample because $a_1\geq 1$ and $b_1 > a_1e$ by hypothesis. Then (a) is equivalent to  the
normal generation  of $L_{1}.$ We will apply~\cite[Theorem 1.3]{rational} for property $N_0$.  Then we just need to show that $-K_{X}\cdot L_{1}=(2C_{0}+(2+e)f)\cdot(a_{1}%
C_{0}+b_{1}f)\geq3.$ This amounts to showing
$-a_{1}e+2a_{1}+2b_{1}\geq 3,$ which is true by hypothesis. This proves (a).
\smallskip

To show the lemma when (b) is satisfied we use~\cite[p. 41, Theorem 2]{Mum}. Then it would be enough to
show  that $H^{1}(L_{2} - L_{1})=H^{1}((a_{2}-a_{1})C_{0}%
+(b_{2}-b_{1})f)$ and $H^{2}(L_{2} -2L_{1})=H^{2}((a_{2}-2a_{1})C_{0}+(b_{2}-2b_{1})f)$ both
vanish. Let $\pi$ $:W\longrightarrow  \mathbf P^1$ be the projection
from $W$ to $\mathbf P^1$. Note that by hypothesis $a_{2}-a_{1} \geq 0$, \ and
$a_{2}-2a_{1}\geq$ $-1$, which implies $R^{1}\pi_* (L_{2}-L_{1})=0$ and $R^{1}\pi_* (L_{2}-2L_{1})=0$. So $H^{j}(L_{2}-jL_{1}%
)=H^{j}(\pi_* (L_{2}-jL_{1}))$ for $j=1,2.$ Since $a_{2}-a_{1}\geq0$, then
 $\pi_* (L_{2}-L_{1})$ splits as the sum of $a_2-a_1+1$ line bundles as follows
\begin{equation}\label{split.pi}
 \pi_* (L_{2}-L_{1})=\mathcal{O}_{\mathbf{P}^{1}}(b_2-b_1) \oplus \cdots \oplus \mathcal{O}_{\mathbf{P}^{1}}(b_2-b_1-(a_2-a_1)e).
\end{equation}
Since by hypothesis $b_{2}-b_{1}\geq (a_2-a_1)e-1$, $H^1\pi_* (L_{2}-L_{1})$, and hence $H^{1}(L_{2}%
-L_{1})$, vanish.  Since a curve has no second cohomology, $H^{2}(L_{2}%
-2L_{1})$ also vanishes. This proves the lemma if (b) holds. 

\smallskip

Now assume (c) holds. We will use Observation~\ref{horace}
repeatedly and apply~\cite[p. 41, Theorem 2]{Mum}. Note that, even though~\cite[p. 41, Theorem 2]{Mum} is stated for ample and base--point--free line bundles in~\cite{Mum}, it is still true  if the line bundle is just base--point--free. 
By  Observation~\ref{horace} it suffices to prove that the
multiplication map
\begin{equation*}
H^{0}(a_1C_0+b_1'f)\otimes H^{0}(f)\overset{\mu'%
}{\longrightarrow }H^{0}(a_1C_0+b_1'f)
\end{equation*}
is surjective for all $b_1' \geq b_1$. Now, since $f$ is base--point--free, ~\cite[p. 41, Theorem 2]{Mum} tells that to prove that $\mu'$ is surjective, it is enough to show $H^{1}(a_1C_0+(b_1'-1)f)=0$
and $H^{2}(a_{1}C_{0}+(b_{1}'-2)f)=0.$ Since $a_{1}>0$, $R^{1}\pi_*(a_1C_0+(b_1'-1)f)$ 
and $R^{1}\pi_*(a_1C_0+(b_1'-2)f)$  both vanish and hence it is enough to prove the vanishing of $H^{1}%
(\pi_*(a_1C_0+(b_1'-1)f))$ and  $H^{2}%
(\pi_*(a_1C_0+(b_1'-2)f))$. The latter cohomology vanishes because a curve has no second cohomology, and the former vanishes arguing as in~\eqref{split.pi} because  
$b_{1}'-1-a_{1}e \geq b_{1}-1-a_{1}e \geq-1$. This settles (c). 

\smallskip

If (d) holds, the lemma follows directly from~\cite[p. 41, Theorem 2]{Mum} arguing as for (b) or (c), having in account that $R^{1}\pi_*((b_2-b_1)f)$ and $R^{1}\pi_*(-C_0+(b_2-2b_1)f)$ vanish and that $b_2-b_1 \geq -1$ by hypothesis. 

\smallskip

We will now assume that (e) holds. Without loss of generality we may assume $a_{1}>0$ 
since in this case $W=\mathbf F_0$  and by
an automorphism of $\mathbf{F}_{0}$ we can interchange $C_{0}$ and $f$ and
rename $a_{2}$ as $a_{1}$ if need be. 
In view
of Observation~\ref{horace}  it is enough to show  that the multiplication maps
\begin{equation*}
\begin{matrix}
H^{0}(a_1'C_0+b_1'f )\otimes
H^{0}(f)&\overset{\mu_{1}}{\longrightarrow }&H^{0}(a_1'C_0+(b_1'+1)f) & \textrm{ and } \\
H^{0}(a_1'C_0+b_1'f )\otimes
H^{0}(C_0)&\overset{\mu_{2}}{\longrightarrow }&H^{0}((a_1'+1)C_0+b_1'f)
\end{matrix}
\end{equation*}
are surjective for all $a_1' \geq a_1$ and all $b_1' \geq b_1$. The maps $\mu_1$ are surjective because (c) holds in this case. The maps $\mu _2$ are also surjective because (c) also holds if $b_1 > 0$, after applying an automorphism of $\mathbf{F}_{0}$ interchanging  $C_{0}$ and $f$.
Finally, if $b_1=0$, after the automorphism of $\mathbf{F}_{0}$, $\mu_2$ becomes 
\begin{equation*}
 H^{0}(a_1'f )\otimes
H^{0}(f)\overset{\mu_{2}}{\longrightarrow }H^{0}((a_1'+1)f).
\end{equation*}
Now by~\cite[p. 41, Theorem 2]{Mum} it is enough to show that 
$H^{1}((a_1'-1)f)=0$ and $H^{2}((a_1'-2)f)=0$. Since $R^{1}\pi_*sf$ 
vanishes for any $s \in \mathbf Z$ and $a_1'-1 \geq a_1-1 \geq 0$, arguing as for (b), (c) or (d) we conclude the surjectivity of $\mu_2$. 
\end{proof}

Now we turn our attention to the  bicanonical morphism of quadruple Galois canonical covers $X$ of smooth rational normal scrolls when $X$ is regular.

\begin{theorem}\label{types34}
Let $\varphi:X\longrightarrow  W$ be a quadruple Galois canonical cover of  Type 3 or 4. Then $2K_{X}$ is
very ample and $|2K_{X}|$ embeds $X$ as a projectively normal variety.
\end{theorem}

\begin{proof} 
Recall that $K_X=\varphi^*(C_0+mf)$ (see Theorem~\ref{quad12}). As in the proof of Theorem~\ref{types12}, since $K_X$ is ample, the normal generation of $2K_X$ is equivalent to the surjectivity of
\begin{equation}\label{projnorm34}
 H^0(2K_X) \otimes H^0(2nK_X) \longrightarrow H^0((2n+2)K_X), 
\end{equation}
for all $n \geq 1$. 
We prove so in two steps. 

\smallskip

\noindent\emph{Step 1.} The 
first step is to show that
\begin{equation*}
H^{0}(2K_{X})\otimes H^{0}(2K_{X})\overset{\beta}{\longrightarrow }H^{0}%
(4K_{X}) %
\end{equation*}
surjects and we argue as in Step 1 of the proof of Theorem~\ref{types12}. We
know (see the structure of $\varphi$ as described in Theorem~\ref{quad12}) that, for both Types 3 and 4, 
\begin{equation*}
\varphi_* \mathcal{O}_{X}=\mathcal{O}_{W}\oplus\mathcal{O}_{W}(-C_{0}%
-(m+1)f)\oplus\mathcal{O}_{W}(-2C_{0}-(e+1)f)\oplus\mathcal{O}_{W}%
(-3C_{0}-(m+e+2)f), 
\end{equation*}
hence we have 
\begin{eqnarray*} 
H^{0}(\varphi_*2K_{X})&=&H^{0}(2C_{0}+2mf)\oplus H^{0}%
(C_{0}+(m-1)f)\oplus H^{0}((2m-e-1)f) \textrm{ \ \ and }\\
H^{0}(\varphi_*4K_{X})&=&H^{0}%
(4C_{0}+4mf)\oplus H^{0}(3C_{0}+(3m-1)f)\oplus H^{0}(2C_{0}+(4m-e-1)f)\oplus \\
&&H^{0}(C_{0}+(3m-e-2)f). 
\end{eqnarray*} 
Then Lemma~\ref{surject.algstruct}, 1) tells that in order to prove the surjectivity of $\beta$ it is enough to show the surjectivity of the
following multiplication maps on $W$:%
\begin{eqnarray*}
 H^{0}(2C_{0}+2mf)\otimes H^{0}(2C_{0}+2mf)&\overset{\beta_{1}}{\longrightarrow 
}&H^{0}(4C_{0}+4mf)\\
 H^{0}(2C_{0}+2mf)\otimes H^{0}(C_{0}+(m-1)f)&\overset{\beta_{2}%
}{\longrightarrow }&H^{0}(3C_{0}+(3m-1)f)\\
  H^{0}(2C_{0}+2mf)\otimes H^{0}((2m-e-1)f)&\overset{\beta_{3}}{\longrightarrow 
}&H^{0}(2C_{0}+(4m-e-1)f)\\
  H^{0}(C_{0}+(m-1)f)\otimes H^{0}((2m-e-1)f)&\overset{\beta_{4}}%
{\longrightarrow }&H^{0}(C_{0}+(3m-e-2)f).
\end{eqnarray*}
To prove the surjectivity of $\beta_1, \beta_2, \beta_3$ and $\beta_4$ we use Lemma~\ref{surjonHirz}. Recall that $m \geq e+1$ and $0 \leq e \leq 2$. Then  
multiplication map $\beta_{1}$ is surjective by Lemma~\ref{surjonHirz}, (a), $\beta_{2}$ by
Lemma~\ref{surjonHirz}, (b), and $\beta_{3}$  and $\beta_{4}$ by Lemma~\ref{surjonHirz}, (c). 

\smallskip

\noindent\emph{Step 2.} To complete the proof, we need to show that the
multiplication map~\eqref{projnorm34}
surjects for all
$n \geq 2$. This follows from the same argument used for Step 2 of the proof of Theorem~\ref{types12}.
\end{proof}

Now we focus on the cases in which $X$ is irregular. Under this hypothesis, in the next theorem we find out when $2K_X$ is normally generated:

\begin{theorem}\label{types5262} Let $\varphi:X \longrightarrow  W$ be a quadruple Galois canonical cover of Type 5.2 or 6.2. Then $2K_{X}$ is
very ample and $|2K_{X}|$ embeds $X$ as a projectively normal variety.
\end{theorem}

\begin{proof}
As we observed in the proof of Theorems~\ref{types12} and~\ref{types34}, the normal generation of $2K_X$ is equivalent to the surjectivity of 
\begin{equation}\label{projnorm5262}
H^0(2K_X) \otimes H^0(2nK_X) \longrightarrow H^0((2n+2)K_X), 
\end{equation}
for all $n \geq 1$. We prove the surjectivity of~\eqref{projnorm5262} in two steps: 

\smallskip

\noindent \emph{Step 1}: 
We show that
\begin{equation*}
H^{0}(2K_{X})\otimes H^{0}(2K_{X})\overset{\beta}{\longrightarrow }H^{0}%
(4K_{X}) %
\end{equation*}
surjects. To see this we argue as in Step 1 of the proofs of Theorems~\ref{types12} and~\ref{types34}. 
Recall that $K_X=\varphi^*(C_0+mf)$ (see Theorem~\ref{quad12}). We
know (see the structure of $\varphi$ as described in Theorem~\ref{quad12})
\begin{equation*}
\varphi_* \mathcal{O}_{X}=\mathcal{O}_{W}\oplus\mathcal{O}_{W}(-C_{0}%
-(m+2)f)\oplus\mathcal{O}_{W}(-2C_{0})\oplus\mathcal{O}_{W}(-3C_{0}-(m+2)f), 
\end{equation*}
hence we have 
\begin{eqnarray*} 
H^{0}(\varphi_*2K_{X})&=&H^{0}(2C_{0}+2mf)\oplus H^{0}%
(C_{0}+(m-2)f)\oplus H^{0}((2mf) \textrm{ \ \ and }\\
H^{0}(\varphi_*4K_{X})&=&H^{0}%
(4C_{0}+4mf)\oplus H^{0}(3C_{0}+(3m-2)f)\oplus H^{0}(2C_{0}+4mf)\oplus \\
&&H^{0}(C_{0}+(3m-2)f). 
\end{eqnarray*} 
Then Lemma~\ref{surject.algstruct}, 1) tells that, in order to prove the surjectivity of $\beta$, it is enough to show the surjectivity of the
following multiplication maps on $W$:%
\begin{eqnarray*}
H^{0}(2C_{0}+2mf)\otimes H^{0}(2C_{0}+2mf)&\overset{\beta_{1}}{\longrightarrow}
&H^{0}(4C_{0}+4mf)\\
H^{0}(2C_{0}+2mf)\otimes H^{0}(C_{0}+(m-2)f)&\overset{\beta_{2}%
}{\longrightarrow}&H^{0}(3C_{0}+(3m-2)f)\\
H^{0}(2C_{0}+2mf)\otimes H^{0}(2mf)&\overset{\beta_{3}}{\longrightarrow}%
&H^{0}(2C_{0}+4mf)\\
H^{0}(C_{0}+(m-2)f)\otimes H^{0}(2mf)&\overset{\beta_{4}}{\longrightarrow}%
&H^{0}(C_{0}+(3m-2)f)
\end{eqnarray*}

To prove the surjectivity of $\beta_1, \beta_2, \beta_3$ and $\beta_4$ we use Lemma~\ref{surjonHirz}. Recall that $m \geq 2$ and $e=0$. Then  $\beta_1, \beta_2, \beta_3$ and $\beta_4$ surject by Lemma~\ref{surjonHirz}, (e). 
\smallskip

\noindent\emph{Step 2}: To complete the proof, we need to show that the
multiplication map~\eqref{projnorm5262}
surjects for all
$n \geq 2$. This follows form the same argument used for Step 2 of the proof of Theorem~\ref{types12}.
\end{proof}

To end the study of whether $2K_X$ is normally generated or not when $W$ is smooth, we look at the remaining types in which $X$ is irregular, i.e, Types 5.1, 6.1, 7, 8. As the following remark shows, for none of these types is $2K_X$ normally generated:

\begin{remark}\label{types516178}
 Let $\varphi:X\longrightarrow  W$ be a quadruple Galois canonical cover of  Type 5.1, 6.1, 7 or 8. Then $|2K_{X}|$ does not embed $X$ as a projectively normal variety.
\end{remark}

\begin{proof}
If $|2K_{X}|$ embedded $X$ as a projectively normal variety, in particular it would embed $X$ as a quadratically normal variety, so the multiplication map
\begin{equation*}
H^{0}(2K_{X})\otimes H^{0}(2K_{X})\overset{\beta}{\longrightarrow }H^{0}%
(4K_{X}) %
\end{equation*}
should surject. Recall that $K_X=\varphi^*(C_0+mf)$ (see Theorem~\ref{quad12}). We are going to apply Lemma~\ref{surject.algstruct}, 2). If $X$ is of type 5.1 or 6.1, then the line bundle $L_2^*$ of~\eqref{split} is $\mathcal O_W(-C_0-3f)$ (see Theorem~\ref{quad12}), so $2C_0+2mf-L_2=C_0-f$ (recall that for Types 5.1 and 6.1, $m=1$). If $X$ is of type 7, the line bundle $L_2^*$ of~\eqref{split} is $\mathcal O_W(-3C_0-f)$ (see Theorem~\ref{quad12}), so $2C_0+2mf-L_2=-C_0+(2m-1)f$. Finally, if $X$ is of type 8, the line bundle $L_2^*$ of~\eqref{split} is $\mathcal O_W(-3C_0)$ (see Theorem~\ref{quad12}), so $2C_0+2mf-L_2=-C_0+2mf$. In any case,  $H^0(2C_0+2mf-L_2)=0$ for Types 5.1, 6.1, 7 and 8. On the other hand, $4C_0+4mf-L_2$ is $3C_0+f$, $C_0+(4m-1)f$ and $C_0+4mf$  for Types 5.1 and 6.1, Type 7 and Type 8 respectively, so in all cases $H^0(4C_0+4mf-L_2) \neq 0$. Then Lemma~\ref{surject.algstruct}, 2) implies that $\beta$ does not surject. 
\end{proof}

\section{The bicanonical morphism when $X$ is irregular}\label{irreg.section}

We have seen in the previous section that if $W$ is smooth and $X$ is regular, then $\varphi_2$ embeds $X$ as a projectively normal variety; in particular, $2K_X$ is very ample. We also saw that if $X$ is irregular, then, in some cases $\varphi_2$ embeds $X$ as a projectively normal variety (Types 5.2 and 6.2) and in others (Types 5.1, 6.1, 7 and 8) does not (see Theorem~\ref{types5262} and Remark~\ref{types516178}). For the latter cases, we will prove in this section that $\varphi_2$ is not even an embedding, we will compute its degree and we will describe $\varphi_2(X)$. Thus in this section $X$ will be irregular and, by Theorem~\ref{quad12}, $W=\mathbf F_0$.

\begin{theorem}\label{Type5-7}
Let $\varphi: X \longrightarrow W$ be a  quadruple Galois canonical cover of Type 5.1, 6.1 or 7. Recall that $K_X=\varphi^*H$ ($H=C_0+mf$; $m =1$ for Types 5.1 and 6.1). Then 
\begin{enumerate}
\item the morphism $\varphi_2$ is a $2:1$; 
\item the image of $\varphi_2$ is a ruled surface of irregularity $m$. More precisely, $\varphi_2(X)= G \times D$, where $G=\mathbf P^1$ and $D$ is a smooth curve of genus $m$ ($m=1$ for Types 5.1 and 6.1), embedded by $|2D+4mG|$ ($G$ is embedded as a conic and $D$ is embedded as a projectively normal curve of degree $4m$). 
\end{enumerate}
\end{theorem}

\begin{proof}
First we deal with $\varphi$ of Types 6.1 and 7. 
Recall that if $\varphi$ is of Type 6.1, then it is the fiber product of two double covers of $W$ branched along $D_1\sim 4C_0$ and $D_2 \sim  2C_0+6f$ 
respectively and if $\varphi$ is of Type 7, then it is the fiber product of two double covers of $W$ branched along $D_1 \sim (2m+2)f$ and $D_2\sim 6C_0+2f$ respectively. Then $\varphi$ factors as $\varphi=p_1 \circ p_2$, where $p_1: X' \longrightarrow W$ is the double cover of $W$ branched along $D_1$ and $p_2:X \longrightarrow X'$ is the  double cover of $X'$ branched along $p_1^*D_2$. 
Thus $\varphi_*\mathcal O_X$ splits as
\begin{equation}\label{split.6.1}
\varphi_*\mathcal O_X=\mathcal O_W \oplus \mathcal O_W(-2C_0)  \oplus \mathcal O_W(-C_0-3f)   \oplus \mathcal O_W(-3C_0-3f) 
\end{equation}
if $\varphi$ is of Type 6.1 and 
\begin{equation}\label{split.dc6.1}
{p_1}_*\mathcal O_{X'} = \mathcal O_W \oplus \mathcal O_W(-2C_0). 
\end{equation}
is the subalgebra of $\varphi_*\mathcal O_X$ that corresponds to $p_1$ in this case. 
Likewise, if $\varphi$ is of Type 7, 
\begin{equation}\label{split.7}
\varphi_*\mathcal O_X=\mathcal O_W \oplus \mathcal O_W(-(m+1)f) \oplus \mathcal O_W(-3C_0-f)  \oplus \mathcal O_W(-3C_0-(m+2)f).  
\end{equation}
and 
the subalgebra of $\varphi_*\mathcal O_X$ that corresponds to $p_1$ is 
\begin{equation}\label{split.dc7}
{p_1}_*\mathcal O_{X'} = \mathcal O_W \oplus \mathcal O_W(-(m+1)f). 
\end{equation}
Since $K_X=\varphi^*H$, it follows from the projection formula and from~\eqref{split.6.1},~\eqref{split.dc6.1},~\eqref{split.7} and~\eqref{split.dc7} that the global sections of $2K_X$ can be identified with the global sections of ${p_1}_*\mathcal O_{X'} \otimes 2H$, so $\varphi_2$ factors through $X'$. More precisely, $\varphi_2=\varphi_2' \circ p_2$, where $\varphi_2'$ is induced by the complete linear series of $L=p_1^*(2H)$ ($m=1$ for Type 6.1). 

\smallskip

\noindent
Now let us study $X'$. The structure of $D_1$ implies that $X'$ is the product of $G=\mathbf P^1$ and a smooth curve $D$. For Type 6.1, $D$ is the pullback of $f$ and has genus $1$. For Type 7,  $D$ is the pullback of $C_0$ and has genus $m$. Using the projection formula and the Leray spectral sequence, it follows easily that $H^1(L-p^*_1C_0)=H^1(L-p_1^*f)=0$, hence $|L|$ restricts to a complete linear series both on $D$ and on $G$. For Type 6.1, the restriction of $L$ to $G$ has degree $2$ and the restriction of $L$ to $D$ has degree $4$, so $\varphi_2'$ embeds $X'$ by $|2D+4G|$ and its image is the product of a smooth conic and a smooth, projectively normal elliptic curve of degree $4$. In particular $\varphi_2$ is $2:1$. 

\smallskip

\noindent Now, for Type 7, the restriction of $L$ to $G$ has degree $2$ as before, so $|L|$ maps $G$ onto a smooth conic. On the other hand the restriction of $L$ to $D$ has degree $4m$, so $|L|$ embeds $D$ as projectively normal curve of degree $4m$. Summarizing, $\varphi_2$ is $2:1$, $\varphi_2'$ is the embedding of $X'$ by the complete linear series $|2D+4mG|$ and the image of $\varphi'_2$ (which is the same as the image of $\varphi_2$) is  the product of a smooth conic and a projectively normal curve of genus $m$ and degree $4m$.

\smallskip
\noindent 
Now we deal with Type 5.1. Recall (see Theorem~\ref{quad12}) that $\varphi$ factors through a double cover $p_1:X' \longrightarrow W$, branched along $D_2 \sim 4C_0$. From the structure of $\varphi$ described in Theorem~\ref{quad12} it follows also that $\varphi_*\mathcal O_X$ splits as~\eqref{split.6.1} and that the subalgebra of $\varphi_*\mathcal O_X$ corresponding to $p_1$ is like the one described in~\eqref{split.dc6.1}.
Then, arguing as for Type 6.1 we see that  $\varphi_2$ factors through $p_1$. Note now that the branch divisor of $p_1:X' \longrightarrow W$ is linearly equivalent to $4C_0$ for both Type 5.1 and 6.1. Thus the rest of the argument for Type 5.1 is the same as for Type 6.1. 
\end{proof}

\begin{theorem}\label{Type8}
Let $\varphi:X \longrightarrow W$ be a  quadruple Galois canonical cover of Type 8. Then 
\begin{enumerate}
\item the surface $X$ is the product of a smooth curve of genus $2$ 
 and a smooth curve of genus $m+1$; 
\item[(2.1)] if $m=1$ (i.e, $\varphi:X \longrightarrow W$ is of Type 8.1), then the morphism $\varphi_2$ is $4:1$ and its image is $W$;
 \item[(2.2)] if $m>1$ (i.e, $\varphi:X \longrightarrow W$ is of Type 8.2), then $\varphi_2$ is $2:1$; 
\item[(3.1)] if $m=1$, then the image of $\varphi_2$ is $\mathbf F_0$, embedded by $|2C_0+2f|$;  
\item[(3.2)] if $m>1$, then the image of $\varphi_2$ is the product of a smooth, projectively normal curve of genus $m+1$ and degree $4m$ and a smooth conic. 
\end{enumerate}
\end{theorem}

\begin{proof}
Recall that $\varphi$ is the fiber product of two double covers of $W$, branched along $D_1\sim (2m+4)f$ and $D_2 \sim 6C_0$ respectively. This implies that $X=G \times D$, where $G$ is a smooth curve of genus $2$ (the reduced structure of the pullback of a fiber $f$ by $\varphi$ at the branch locus of $\varphi$)
 and $D$ is a smooth curve of genus $m+1$ (the reduced structure of the pullback of $C_0$ by $\varphi$ at the branch locus of $\varphi$). We will show that $|2K_X|$ restricts to a complete linear series both on $G$ and on $D$. For this it suffices to show that $H^1(2K_X-G)=H^1(2K_X-D)=0$. To see $H^1(2K_X-G)=0$ we use Kodaira vanishing, noting that $G$ is numerically equivalent to $\frac 1 2 \varphi^*f$, so $K_X-G$ is numerically equivalent to $\varphi^*(C_0+(m-\frac 1 2 f))$, hence ample. To see $H^1(2K_X-D)=0$ we argue similarly, noting that $D$ is numerically equivalent to $\frac 1 2 \varphi^*C_0$. Then $\varphi_2(X)$ is the product of $\varphi_2(G)$ and $\varphi_2(D)$.

  \smallskip
\noindent 

For general $f$, $\varphi^{-1}f$ is the disjoint union of two curves, each algebraically equivalent to $G$. Thus $G^2=0$. Likewise,  for general $C_0$,  $\varphi^{-1}C_0$ is the disjoint union of two curves, each algebraically equivalent to $D$, so $D^2=0$. Then $2K_X|_{_G}=2K_G$ and $2K_X|_{_D}=2K_D$. Since the genus of $G$ is $2$, then $\varphi_2|_G$ is $2:1$ and its image is a smooth conic. Since the genus of $D$ is $m+1$, if $m >1$, then $\varphi_2|_D$ is an embedding, and if $m=1$, then $\varphi_2|_D$ is $2:1$ and its image is a smooth conic. Thus if $m=1$, $\varphi_2$ is $4:1$ and its image is $\mathbf F_0$, embedded by $|2C_0+2f|$ in $\mathbf P^8$. Finally, if $m >1$, then $\varphi_2$ is $2:1$ and its image is ruled surface of irregularity $m+1$; precisely, the product of a smooth, projectively normal curve of genus $m+1$ and degree $4m$ and a smooth conic. This completes the proof. 
\end{proof}

\begin{remark}\label{Type8.remark} {\rm To prove Theorem~\ref{Type8} 
we could have also argued as in the proof of Theorem~\ref{Type5-7}, as we outline now. Let us call as $p_1:X' \longrightarrow W$  the double cover of $W$ branched along $D_1$. 
Note that 
\begin{equation}\label{split.8}
\varphi_*\mathcal O_X=\mathcal O_W \oplus \mathcal O_W(-(m+2)f) \oplus \mathcal O_W(-3C_0)  \oplus \mathcal O_W(-3C_0-(m+2)f).  
\end{equation}
and that 
\begin{equation}\label{split.dc8}
{p_1}_*\mathcal O_{X'}=\mathcal O_W \oplus \mathcal O_W(-(m+2)f)
 \end{equation}
is the subalgebra that corresponds to $p_1$. Then the projection formula implies that the global sections of $2K_X$ can be identified with the global sections of $2H \otimes (\mathcal O_W \oplus \mathcal O_W(-(m+2)f))$. Moreover, for Type 8.1, the global sections of $2K_X$ can be identified with the global sections of $2H$, since in this case $m=1$. This implies that for Type 8.1, $\varphi_2$ factors through $\varphi$ (and it is therefore $4:1$) and for Type 8.2, $\varphi_2$ factors through $p_1$.  Now the argument would follow the same lines as the proof of Theorem~\ref{Type5-7}. {\hfill $\square$ \medskip}}
\end{remark}

Theorems~\ref{Type5-7} and~\ref{Type8} say 
in particular that if $X$ is of Type 5.1, 6.1, 7 or 8, $\varphi_2$ is not birational. There is a good reason for this to happen, namely, the fact that in these cases $X$ has a 
pencil of curves of genus $2$. In the next remarks we observe that, indeed, if $X$ has a pencil of curves of genus $2$, $\varphi_2$ cannot be birational and we explicitly describe the pencil of curves of genus $2$ if $X$ is of Type 5.1, 6.1, 7 or 8. 

\begin{remark} If $X$ has a pencil of genus $2$ curves, then this pencil is base--point--free and $\varphi_2$ is not birational. 
\end{remark}

\begin{proof}
The pencil of genus $2$ curves is base--point--free because of our hypothesis of $K_X$ being ample and base--point--free. 
Then, if $C$ is a general member of the pencil, $|2K_X|$ restricts to a linear base--point--free subseries of $|K_C|$, hence to the complete canonical linear series of $C$. Then $\varphi_2$ maps $C$ two--to--one onto its image. 
\end{proof}

\begin{remark} We exhibit explicitly a pencil of genus $2$ curves on $X$ of Type 5.1, 6.1, 7 and 8. In particular
\begin{enumerate}
\item if $X$ is of Type 5.1 or Type 6.1, $X$ possesses an elliptic pencil of curves of genus $2$;
\item if $X$ is of Type 7, $X$ possesses a genus $m$ pencil of curves of genus $2$;
\item if $X$ is of Type 8, $X$ possesses a genus $m+1$ pencil of curves of genus $2$;
\end{enumerate}
\end{remark}

\begin{proof} 
If $X$ is of Type 5.1, 6.1, 7 we saw in the proof of Theorem~\ref{Type5-7}
that $\varphi=p_1 \circ p_2$, where $p_1: X' \longrightarrow W$ and $p_2: X \longrightarrow X'$ are double covers, $X'=D \times \mathbf P^1$  and $D$ is a smooth curve of genus $1$ if $X$ is of Type 5.1 or 6.1, and a smooth curve of genus $m$ if $X$ is of Type 7. Then $p_2$ is branched along a divisor meeting a general fiber of $X'$ at $6$ six distinct points (see Theorem~\ref{quad12}),
so $X$ is a fibration over $D$ whose general fiber is a smooth curve of genus $2$.

\smallskip
\noindent 
If $X$ is of Type 8, we proved in  Theorem~\ref{Type8} that $X$ is the product of a smooth curve of genus $2$ and a smooth curve of genus $m+1$, so obviously $X$ possesses a genus $m+1$ base--point--free pencil of curves of genus $2$. 

\smallskip
\noindent 
Note
that the fibration above comes from the Stein factorization of $g$, where $g$ is the composition of $\varphi$ followed by the projection of $W$ onto $f$ if $X$ is of Type 5.1 or 6.1, and by the projection of $W$ onto $C_0$ if $X$ is of Type 7 or 8. \end{proof}

\medskip

\begin{remark} {\rm The existence of pencils of genus $2$ curves if $X$ is of Type 5.1, 6.1, 7 or 8 follows also indirectly from more general results on the classification of the bicanonical maps of surfaces of general type, stated in ~\cite{Xiao}, having in account what the degree and how the image of $\varphi_2$ are (see Theorems~\ref{Type5-7} and~\ref{Type8}). 

\smallskip
\noindent Indeed; if  $X$ is of Type 5.1, 6.1, 7 or 8.2, then Theorems~\ref{Type5-7} and~\ref{Type8}
say that $\varphi_2$ has degree $2$. If $X$ did not have a pencil of genus $2$ curves, \cite[Theorem 2]{Xiao} would imply that either $q(X)=0$ or $\varphi_2(X)$ is rational. However Theorem~\ref{quad12}
says that $q(X)>0$. On the other hand, Theorems~\ref{Type5-7} and~\ref{Type8}
say that the image of $\varphi_2$ is a non--rational ruled surface, so $X$ should have a pencil of genus $2$ curves.

\smallskip
\noindent If $X$ is of Type 8.1, Theorem~\ref{quad12}
implies that $X$ does not fit in the list of exceptions of~\cite[Theorem 1]{Xiao}. On the other hand Theorem~\ref{Type8}
says that $\varphi_2$ has degree $4$, hence~\cite[Theorem 1]{Xiao} implies that $X$ should have a pencil of genus $2$ curves.}
\end{remark}

\medskip

\noindent In Theorems~\ref{Type5-7} and~\ref{Type8} we have described  the image of $\varphi_2$ for $X$ of each of Types 5.1, 6.1, 7 and 8. We settle now the question of whether these images are projectively normal varieties: 

\begin{cor}\label{image.N0}
Let $\varphi:X \longrightarrow W$ be a quadruple Galois canonical cover of Type 5.1, 6.1, 7 or 8. The image of $X$ by $\varphi_2$ is a projectively normal variety. 
\end{cor}

\begin{proof}
If $\varphi$ is of Type 8.1, $\varphi_2(X)$ is $\mathbf F_0$ embedded by $|2C_0+2f|$, which is well known to be projectively normal (it does even satisfy property $N_5$; see~\cite[Theorem 1.3]{rational}). For the remaining types recall that $\varphi_2(X)=\varphi_2'(X')$, where $p_1:X' \longrightarrow W$ is a double cover branched along a divisor linearly equivalent to $4C_0$ in Types 5.1 and 6.1, to $(2m+2)f$, in Type 7 and to $(2m+4)f$ in Type 8.2. Recall also that $\varphi'_2$ is induced by $|L|$, where $L=p_1^*(2H)$. Since $L$ is ample, to prove that $|L|$ embeds $X'$ as a projectively normal variety it suffices to show that the 
multiplication maps 
\begin{equation*}
H^0(L) \otimes H^0(rL) \longrightarrow H^0((r+1)L)
\end{equation*}
are surjective for all $r \geq 1$. 
Arguing as in Section~\ref{Normal.generation.2K_X} 
it follows from~\eqref{split.dc6.1},~\eqref{split.dc7} and~\eqref{split.dc8} that it suffices to check that the following multiplication maps on $W$ surject:
\begin{equation}\label{surj.5.1.6.1}
\begin{matrix}
H^0(2C_0+2f) \otimes H^0(2rC_0+2rf) & \longrightarrow &  H^0((2r+2)C_0+(2r+2)f)  
& \text{and} \\
\\
H^0(2f) \otimes H^0(2rC_0+2rf)  & \longrightarrow & H^0((2rC_0+(2r+2)f), 
\end{matrix}
\end{equation}
for Types 5.1 and 6.1; 
\begin{equation}\label{surj.7}
\begin{matrix}
H^0(2C_0+2mf) \otimes H^0(2rC_0+2rmf) & \longrightarrow &  H^0((2r+2)C_0+(2r+2)mf) & \text{and}\\
\\
H^0(2C_0+(m-1)f) \otimes H^0(2rC_0+2rf)  & \longrightarrow & H^0((2rC_0+(2r+m-1)f),
\end{matrix}
\end{equation}
for Type 7; and
\begin{equation}\label{surj.8}
\begin{matrix}
H^0(2C_0+2mf) \otimes H^0(2rC_0+2rmf) & \longrightarrow &  H^0((2r+2)C_0+(2r+2)mf) & \text{and}\\
\\
H^0(2C_0+(m-2)f) \otimes H^0(2rC_0+2rf)  & \longrightarrow & H^0((2rC_0+(2r+m-2)f),
\end{matrix}
\end{equation}
($m>1$), for Type 8. 
By Observation~\ref{horace}, 
to prove the surjectivity of the maps~\eqref{surj.5.1.6.1},~\eqref{surj.7} and~\eqref{surj.8} it suffices to 
see that the multiplication maps on $W$ 
\begin{equation}
\begin{matrix}
H^0(s_1C_0+s_2f) \otimes H^0(C_0) & \longrightarrow &  H^0((s_1+1)C_0+s_2f) & \text{and} \\
\\
H^0(s_1C_0+s_2f) \otimes H^0(f) & \longrightarrow &  H^0(s_1C_0+(s_2+1)f)
\end{matrix}
\end{equation}
surject for all $s_1,s_2 \geq 2$. 
This follows from~\cite[p. 41, Theorem 2]{Mum}.
\end{proof}

\section{The bicanonical morphism when $W$ is singular}\label{sing.section}

In this section we prove that if $W$ is singular (and hence, see~\cite{GPsing}, $W=S(0,2)$) then $2K_X$ is not very ample. We also see that, in constrast with the case when $W$ is smooth (in which either $\varphi_2$ is an embedding with projectively normal image or $\varphi_2$ has degree bigger than $1$) if $W$ is singular,  then $\varphi_2$ is birational but not an embedding. 
Precisely, we have this 

\begin{theorem}\label{sing}
Let ${\varphi}:X  \longrightarrow W$ be a quadruple Galois canonical cover 
of $W=S(0,2)$ and let $w$ the vertex of $W$.

Then $\varphi_2$ is a birational morphism but $2K_X$ is not very ample because $\varphi_2$ fails to be an embedding at $\varphi^{-1}\{w\}$. 
Precisely, 
\begin{enumerate}
\item if ${\varphi}$ is of Type 9, 10 or 12, then $|2K_X|$ does not separate directions at the unique point $x \in \varphi^{-1}\{w\}$;  
\item if ${\varphi}$ is of Type 11 then $|2K_X|$ does not separate the two points $x_1$ and $x_2$ of $\varphi^{-1}\{w\}$, although ${\varphi_2}$ is locally an embedding at both of them. 
\end{enumerate}
Moreover, $\varphi_2$ is an embedding on the set of smooth points of $X-\varphi^{-1}\{w\}$ and, if $X$ has the mildest possible singularities (see~\cite[Corollary 5.1 and Propositions 5.2, 5.3 and 5.4]{GPsing}), then $\varphi_2$ is an embedding outside $\varphi^{-1}\{w\}$. 
\end{theorem}

Before we prove Theorem~\ref{sing} we state a couple of lemmas. 

\begin{lemma}\label{blowup}
Let $S$ be an irreducible, normal surface with only rational singularities, let $\pi: \widetilde S \longrightarrow S$ be a proper birational map and let $\mathfrak m$ be a maximal ideal sheaf on $S$. Then
$\pi_*(\mathfrak m^n\mathcal O_{\widetilde S})=\mathfrak m^n$. 
\end{lemma}

\begin{proof} The lemma follows from Remark c) of Section 5, Theorem 7.1 and Proposition 6.2 of~\cite{Lip}.  

\end{proof}

\begin{lemma}\label{hyper}
Let $S$ be a regular surface of general type whose canonical map is a degree $n$ morphism $\psi:S \longrightarrow W$ onto a surface of minimal degree $W$.
Let $C \in |K_S|$ be a smooth, irreducible curve. 
If $C$ is hyperelliptic, then $n=2$. 
\end{lemma}

\begin{proof}
Since $K_S$ is base--point--free, then $K_S|_{_C}$ is a base--point--free theta--characteristic on $C$. 
Let $g$ be the genus of $C$. 
Thus, if $C$ is hyperelliptic, $K_S|_{_C}$ is $\frac{g-1}{2}$ times the $g^1_2$ of $C$ and $h^0(K_S|_{_C})=\frac{g+1}{2}$. 
On the other hand, since $S$ is regular, $h^0(K_S|_{_C})=p_g-1$. Since $W$ is a surface of minimal degree, $g-1=$deg$K_S|_{_C}=K_S^2=n(p_g-2)$, hence $g-1=n(p_g-2)=n(\frac{g+1}{2}-1)$, so $n=2$. 
\end{proof}

\begin{Proofof} \emph{Theorem~\ref{sing}.} We use the notation of~\ref{CD}. 
We first prove that if $X$ has the mildest possible singularities, then $\varphi_2$ is an embedding outside $\varphi^{-1}\{w\}$. 
Recall that outside $w$, $Y$ and $W$ are isomorphic. Recall also that outside $\varphi^{-1}\{w\}$, $\overline X$ and $X$ are isomorphic and $p$ and $\varphi$ are equal. Let now $x_1$ and $x_2$ be two distinct  points of $X - \varphi^{-1}\{w\}$.  
If $\varphi(x_1) \neq \varphi(x_2)$, then $\varphi_2(x_1) \neq \varphi_2(x_2)$, since $\varphi_2(X)$ can be projected to a $2$--Veronese reembedding of $W$. Thus, let us assume $\varphi(x_1) = \varphi(x_2)$ and let us call $y$ the inverse image of $\varphi(x_1) = \varphi(x_2)$ by $q$. Consider the linear system $|C_0+2f|$ in $Y$.
Since it is base--point--free and big,  
there exists a smooth connected curve in $|C_0+2f|$, 
not meeting $C_0$, passing through $y$ and meeting the branch locus of $p$ in such a way that its pullback $C$ by $p$ is smooth and irreducible. 
Then $\overline q(C)$ is isomorphic to $C$ and belongs to $|K_X|$, is smooth and irreducible and passes through $x_1$ and $x_2$. Since $\varphi$ has degree $4$, $C$ is non hyperelliptic by Lemma~\ref{hyper}. Since $X$ is regular, $|2K_X|$ restricts to the complete canonical series of $C$, thus $|2K_X|$ embeds $C$ and therefore separates $x_1$ and $x_2$. A similar argument proves that in general $\varphi_2$ separates smooth points of $\overline X-\varphi^{-1}\{w\}$. Moreover, this argument can be adapted to show that $\varphi_2$ separates directions in the set of smooth points of $\overline X-\varphi^{-1}\{w\}$.

\medskip

\noindent Now we see what happens at $\varphi^{-1}\{w\}$. First we treat the case in which $\varphi^{-1}\{w\}$ consists of only one point $x$ (Types 9, 10 and 12). Let $\mathcal O_{X,x}$ be the local ring at $x$, let $\mathfrak m_x$ be the maximal ideal of  $x$  
and let $$Z=\mathrm{\ Spec\ }\mathcal O_{X,x}/\mathfrak m_x^2$$
 be the first infinitesimal neighbourhood of $x$. We want to prove that $|2K_X|$ does not separates directions at $x$. This is equivalent to proving that 
$$H^0(2K_X) \longrightarrow H^0(2K_X|_{_Z})$$
 is not surjective.   Since, by Kodaira vanishing $H^1(2K_X)=0$, the latter is equivalent to the non vanishing of $H^1(2K_X \otimes \mathfrak m_x^2)$. 
To study this cohomology group we will use Lemma~\ref{blowup}. Thus 
we will compute $H^1(2K_X \otimes \mathfrak m_x^2)$ by arguing on $\overline X$ and with $L=\overline q^*K_X=p^*(C_0+2f)$. We denote $F=p^{-1}C_0$.

\medskip

\noindent Now we argue for each Type 9, 10 and 12. We start with Type 12. 
In this case (see~\cite[Corollary 5.1]{GPsing}) 
the morphism $\overline X \overset{\overline q}\longrightarrow X$ is the
contraction of the smooth line $F$. 
The line $F$
consists of smooth points of $\overline X$ and an $A_1$
singularity $\overline x$. Recall also that the self--intersection of $F$ is $F^2=-\frac
1 2$ and that $p^*C_0=4F$.
Thus $\overline q$ factorizes as a composition of blowing ups and blowing downs in the following way: Let $g: \widehat X \longrightarrow \overline X$ be the blowing up of $\overline X$ at $\overline x$. Then the exceptional divisor $F_1$ of $g$ has $F_1^2=-2$ and the strict transform of $F$ is a line $F_2$ with $F_2^2=-1$. Then we obtain $X$ from $\widehat X$ by contracting first $F_2$ and then $F_1$, so  $\widehat X$ is obtained from $X$ by performing two consecutive blowing ups, $g_1:X' \longrightarrow X$ and $g_2 : \widehat X \longrightarrow X'$, the first one centered at $x$. Let us call $M=g_2^*(g_1^*K_X)$ and recall that $L=\overline q^*K_X=p^*(C_0+2f)$. Notice also that $M=g^*L$. Then, local computations and Lemma~\ref{blowup} yield 
$$(g_2 \circ g_1)_*\mathcal O_{\widehat X}(-2F_1-2F_2)=\mathfrak m_x^2.$$ 

Then, by projection formula, 
$$(g_2 \circ g_1)_*(2M \otimes \mathcal O_{\widehat X}(-2F_1-2F_2))=2K_X \otimes \mathfrak m_x^2$$
 and by the Leray Spectral Sequence, 
 $$H^1((g_2 \circ g_1)_*(2M \otimes \mathcal O_{\widehat X}(-2F_1-2F_2)))=H^1(2M \otimes \mathcal O_{\widehat X}(-2F_1-2F_2)).$$ 
Local computation shows that
\begin{eqnarray*}
\mathcal O_{\overline X}(-2F)\cdot\mathcal O_{\widehat X}= \mathcal O_{\widehat X}(-F_1-2F_2) \ \mathrm{and \ 
that} \cr
(\mathcal O_{\overline X}(-2F) \otimes m_{\overline x})\cdot \mathcal O_{\widehat X}= \mathcal O_{\widehat X}(-2F_1-2F_2),
\end{eqnarray*} 

where $\mathfrak m_{\overline x}$ is the maximal ideal corresponding to $\overline x$. 

Then Remarks c) and e) of Section 5 and Proposition 6.2 of~\cite{Lip} and Lemma~\ref{blowup} show. 
that 
\begin{eqnarray*}
g_*\mathcal O_{\widehat X}(-F_1-2F_2)=\mathcal O_{\overline X}(-2F) \mathrm{\ and \ } \cr
g_*\mathcal O_{\widehat X}(-2F_1-2F_2)=\mathcal O_{\overline X}(-2F) \otimes \mathfrak m_{\overline x}
\end{eqnarray*}
where $\mathfrak m_{\overline x}$ is the maximal ideal corresponding to $\overline x$. Thus we consider the exact sequence
\begin{equation*}
0 \longrightarrow  2M \otimes \mathcal O_{\widehat X}(-2F_1-2F_2) \longrightarrow 2M \otimes \mathcal O_{\widehat X}(-F_1-2F_2) \longrightarrow 2M \otimes \mathcal O_{F_1}(-F_1-2F_2) \longrightarrow 0
\end{equation*}
and push it down to $\overline X$ to obtain
\begin{equation}\label{seq} 
0 \longrightarrow  g_*(2M \otimes \mathcal O_{\widehat X}(-2F_1-2F_2)) \longrightarrow 2L \otimes \mathcal O_{\overline X}(-2F)\longrightarrow \mathbf k_{\overline x} \longrightarrow 0 , 
\end{equation}
where $\mathbf k_{\overline x}$ is the skycraper sheaf on $\overline x$ of dimension $1$ obtained by restricting $2M \otimes \mathcal O_{F_1}(-F_1-2F_2)$ to $\overline x$. The above sequence is exact because $R^1g_*(2M \otimes \mathcal O_{\widehat X}(-2F_1-2F_2))$ vanishes. 
Again by Leray Spectral Sequence, $H^1(2M \otimes \mathcal O_{\widehat X}(-2F_1-2F_2))=H^1(g_*(2M \otimes \mathcal O_{\widehat X}(-2F_1-2F_2)))$, so we take cohomology in the above sequence~\eqref{seq}. So, if we see that $H^1(2L \otimes \mathcal O_{\overline X}(-2F)) \neq 0$, then
$H^1(2M \otimes \mathcal O_{\widehat X}(-2F_1-2F_2)) \neq 0$ and we are done. 
Then, to see that $H^1(2L \otimes \mathcal O_{\overline X}(-2F)) \neq 0$, we argue like this. 
From~\cite[Corollary 5.1]{GPsing},  $K_{\overline
X}={\overline q}^*K_X + 2F$, hence  
$$H^1(2L \otimes \mathcal O_{\overline X}(-2F))=H^1(K_{\overline
X} + L-4F)=H^1(K_{\overline
X}^* +p^*(2f))
=H^1(p^*(-2f))^*.$$ 
Now by~\cite[Corollary 4.3]{GPsing}, 
\begin{equation}\label{seq.splitting0}
p_*\mathcal O_{\overline X}=\mathcal O_Y \oplus \mathcal O_Y(-2C_0-3f) \oplus \mathcal O_Y(-2C_0-3f) \oplus \mathcal O_Y(-3C_0-6f),
\end{equation}
so 
\begin{equation}\label{seq.splitting}
p_*p^*(\mathcal O_{Y}(-2f))= \mathcal O_Y(-2f) \oplus \mathcal O_Y(-2C_0-5f) \oplus \mathcal O_Y(-2C_0-5f) \oplus \mathcal O_Y(-3C_0-8f)
\end{equation}
Then, by~\eqref{seq.splitting}
and the Leray Spectral Sequence, $h^1(p^*(-2f))=h^1(\mathcal O_X(-2f))=1$.

\medskip
\noindent Now we deal with Type 10. In this case $\overline q$ is the blow--up of $X$ at $x$ and a partial desingularization of $X$ at $x$ (recall that $x$ is a $D_4$ singularity). Recall also that the exceptional divisor of $\overline q$ is a line $F$ and in this case we have, as in Type 12, that $p^*C_0=4F$. The points of $F$ are smooth except $3$ points which are $A_1$ singularities and $F^2=-1/2$. 
The local equation of $X$ at $x$ is $z^2t-t^3-u^2=0$ and a local computation of the blowing up at $x$ shows that $\mathfrak m_x\mathcal O_{\overline X}=\mathcal O_{\overline X}(-2F)$ 
and $\mathfrak m_x^2\mathcal O_{\overline X}=\mathcal O_{\overline X}(-4F)$. 
Now recall that we want to prove the nonvanishing of $H^1(2K_X \otimes \mathfrak m_x^2)$. Recall that $L=\overline q^*K_X=p^*(C_0+2f)$. By Lemma~\ref{blowup}, $\overline q_*(\mathcal O_{\overline X}(-4F))=\mathfrak m_x^2$, hence, by the projection formula and the Leray spectral sequence, 
$$H^1(2K_X \otimes \mathfrak m_x^2)=H^1(p^*(C_0+4f)).$$
Recall also that in Type 10 (see~\cite[Corollary 4.3]{GPsing}),  
\begin{equation}\label{seq.splitting2}
p_*\mathcal O_{\overline X} = \mathcal O_Y \oplus \mathcal O_Y(-C_0-3f) \oplus \mathcal O_Y(-2C_0-3f) \oplus \mathcal O_Y(-3C_0-6f). 
\end{equation}
Thus, by the projection formula, the Leray spectral sequence and Serre duality, 
$h^1(p^*(C_0+4f))=h^1(\mathcal O_Y(-2C_0-2f))=h^1(\mathcal O_Y(-2f))=1$, so $H^1(2K_X \otimes \mathfrak m_x^2) \neq 0$, as wanted. 

\medskip

\noindent Now we deal with Type 9. In this case, $p^*C_0=2F$ and $F^2=-2$, so $x$ is an $A_1$ singularity, and $\overline q$ is the blow--up of $X$ at $x$, which desingularizes $X$.  Then 
$\mathfrak m_x^2\mathcal O_{\overline X}= \mathcal O_{\overline X}(-2F)$, so arguing as in the cases above and using Lemma~\ref{blowup}, the projection formula and the Leray spectral sequence we have that $H^1(2L\otimes\mathcal O_{\overline X}(-2F))=H^1(2K_X\otimes \mathfrak m_x^2)$. Recall that $L=\overline q^*K_X=p^*(C_0+2f)$, so 
$H^1(2L\otimes\mathcal O_{\overline X}(-2F))=H^1(p^*(C_0+4f))$. Since $p_*\mathcal O_{\overline X}$ splits as~\eqref{seq.splitting2} (see~\cite[Corollary 4.3]{GPsing}), we obtain that 
$h^1(\mathcal O_Y(-2f))=1$ like in Type 10.

\medskip 

\noindent 
Finally we study the case in which $\varphi^{-1}\{w\}$ consists of two points, $x_1$ and $x_2$. This is a quadruple Galois canonical cover of Type 11. We prove first that $|2K_X|$ does not separate $x_1$ and $x_2$. Recall  (see~\cite[Corollary 5.1]{GPsing}) that $x_1$ and $x_2$ are smooth points and $\overline q$ is the blowing up of $X$ at $x_1$ and $x_2$, so $\overline q_*\mathcal O_{\overline X}(-F_1-F_2)=\mathfrak m_{x_1} \otimes \mathfrak m_{x_2}$, where 
$F_1$ and $F_2$ are the exceptional divisors of $\overline q$, which are $-1$--lines. Let $f$ be a general fiber of the ruled surface $Y$ and let $\overline f$ be the pullback to $\overline X$ of $f$ by $p$. Then $\overline f$ is a smooth, connected curve
of genus $4$, meeting $F_1$ (respectively $F_2$) at one point $\overline x_1$ (respectively $\overline x_2$) transversally (see the proof of~\cite[Theorem 4.1]{GPsing}).  
Recall also that $\overline q^*K_X=L=p^*(C_0+2f)$ and that the morphism induced on $\overline X$ by $|2L|$ factors through $\varphi_2$. Then, if the restriction of $|2L|$ to $\overline f$ does not separate 
$\overline x_1$ and $\overline x_2$, then $\varphi_2(x_1)=\varphi_2(x_2)$. Note that the degree of $2L|_{\overline f}$ is $8$. Recall also that $K_{\overline X}=p^*(C_0+2f)+F_1+F_2$. Then, by adjunction formula,  $2L|_{\overline f}$ is the canonical of $\overline f$ plus the degree $2$, effective divisor $(F_1+F_2)|_{\overline f}=\overline x_1+\overline x_2$. Therefore the restriction of $|2L|$ to $\overline f$ does not separate 
$\overline x_1$ and $\overline x_2$. 

\smallskip

\noindent Now we 
show that $\varphi_2$ is a local embedding at both $x_1$ and $x_2$. 
Let $Z'=$Spec$(\mathcal O_{X,x_1}/\mathfrak m_{x_1}^2 \oplus \mathcal O_{X,x_2}/\mathfrak m_{x_2}^2)$. 
We will show that the cokernel of the homomorphism
\begin{equation}\label{restrZ}
H^0(2K_X) \longrightarrow H^0(2K_X|_{Z'})
\end{equation}
has dimension $1$. Since $\varphi_2(x_1)=\varphi_2(x_2)$, in that case it would follow that $\varphi_2$ is a local embedding at both $x_1$ and $x_2$, so we would be done. 

Thus, we show now that the cokernel of~\eqref{restrZ} has dimension $1$. Since $H^1(2K_X)=0$, this is equivalent to showing that $h^1(2K_X \otimes \mathfrak m_{x_1}^2 \otimes \mathfrak m_{x_2}^2)=1$. For the latter, since $\overline q_*\mathcal O_{\overline X}(-2F_1-2F_2)=\mathfrak m_{x_1}^2 \otimes \mathfrak m_{x_2}^2$, from the projection formula and the Leray spectral sequence it follows 
\begin{equation*}
H^1(2K_X \otimes \mathfrak m_{x_1}^2 \otimes \mathfrak m_{x_2}^2)=H^1(2L \otimes \mathcal O_{\overline X}(-2F_1-2F_2)).
\end{equation*} 
The latter cohomology group is equal to $H^1(p^*(C_0+4f))$. The splitting of $p_*\mathcal O_{\overline X}$ is~\eqref{seq.splitting0} (see~\cite[Corollary 4.3]{GPsing}), so the projection formula, the Leray spectral sequence and Serre duality yield 
$h^1(p^*(C_0+4f))=h^1(\mathcal O_Y(-2C_0-4f))=h^1(\mathcal O_Y(-2f))=1$, so $h^1(2K_X \otimes \mathfrak m_{x_1}^2 \otimes \mathfrak m_{x_2}^2)=1$.
\end{Proofof}

\section{Ring generators of the canonical ring of quadruple canonical
covers}\label{canring.section}

In this section we study the  generators of the canonical ring of a quadruple Galois canonical cover $X$ of minimal degree. Precisely we will find the degrees of the minimal generators of the canonical ring of $X$, and the number of generators in each degree. 
The first result  of this section is a general result that gives a beautiful formula for the number of generators in degree $2$ of the canonical ring of canonical covers, of arbitrary degree and irregularity, of surfaces of minimal degree. This result recovers part.
of~\cite[Theorem 2.1]{tams}, that gives the degrees and number of generators of the canonical ring 
of a regular surface $S$ of general type with at worst canonical singularities, that is a canonical cover of a surface of minimal degree.
One of the consequences of~\cite[Theorem 2.1]{tams} is that the number of extra generators needed in degree $2$ depends only on the geometric genus $p_g(S)$ of $S$ and on the degree $n$ of $\Psi$; precisely~\cite[Theorem 2.1]{tams} tells that this number is $(n-2)(p_g(S)-2)$. This formula is generalized to canonical covers  of arbitrary irregularity of surfaces of minimal degree in the following:

\begin{theorem}\label{gen.in.deg2}
 Let $S$ be a surface of general type, normal with at worst canonical singularities and such that its canonical bundle is base--point--free. Let $\Psi$ be the canonical morphism of $S$ and let $n$ be the degree of $\Psi$. Assume that the image of $\Psi$ is a surface $Y$ of minimal degree. Then, the degree $2$ part of the canonical ring $R(S)$ of $S$ is generated by the elements of degree $1$ of $R(S)$ and by $(n-2)(p_g(S)-2)-q(S)$ linearly independent elements of degree $2$.
\end{theorem}

\begin{proof}
The part $R_2$ of degree $2$ of $R(S)$ is $H^0(2K_S)$. Then the number of linearly independent elements of degree $2$ which, in addition to the elements of degree $1$, are needed  to generate $R_2$ is the dimension of the cokernel of the multiplication map of global sections on $S$ 
\begin{equation*}
 H^0(K_S) \otimes H^0(K_S) \longrightarrow H^0(2K_S), 
\end{equation*}
which is equal to the dimension of the cokernel of the multiplication map of global sections on $Y$
\begin{equation*}
 H^0(\Psi_*K_S) \otimes H^0(\Psi_*K_S) \overset{\gamma'}\longrightarrow H^0(\Psi_*2K_S). 
\end{equation*}
 We have an exact sequence 
\begin{equation}\label{pushdown.seq}
 0 \longrightarrow \mathcal O_Y \overset{i}\longrightarrow \Psi_*\mathcal O_S \longrightarrow \mathcal F \longrightarrow 0, 
\end{equation}
of sheaves on $Y$, where $\mathcal F$ is simply the cokernel of $i$. Let $H$ be the hyperplane section of $Y$. Then $K_S=\Psi^*H$ and $H^0(K_S)=H^0(\Psi_*K_S)=H^0(H)$, last equality being induced by $i$. On the other hand, 
\begin{equation*}
 H^0(H) \otimes H^0(H) \longrightarrow H^0(2H)
\end{equation*}
surjects, because $Y$ is projectively normal. 
Then the image of $\gamma'$ is $H^0(2H)$, which by~\eqref{pushdown.seq} is a subspace of $H^0(\Psi_*2K_S)$. Therefore the dimension of the cokernel of $\gamma'$ is $h^0(\Psi_*2K_S)-h^0(2H)=h^0(2K_S)-h^0(2H)$. By Riemann--Roch
and the Kawamata--Viehweg vanishing theorem $h^0(2K_S)=K_S^2+p_g(S)-q(S)+1$. Let now $C$ be a smooth curve in the linear system $H$, not meeting the singular locus of $Y$. Then $H^0(2H)$ fits in the following
\begin{equation*}
 0 \longrightarrow H^0(H) \longrightarrow H^0(2H) \longrightarrow H^0(2H|_C) \longrightarrow 0 
\end{equation*}
exact sequence, since $H^1(H)=0$, which follows from the fact that $Y$ is a regular surface.
Since $C=\mathbf P^1$, $h^0(2H|_C)=2H^2+1$, so $h^0(2H)=p_g(S)+2H^2+1$. Now we have that $K_S^2=nH^2$
 and, since $Y$ is a surface of minimal degree, $H^2=p_g(S)-2$. All this yields $h^0(2K_S)-h^0(2H)=(n-2)(p_g(S)-2)-q(S)$
\end{proof}

The result \cite[Theorem 2.1]{tams} tells among other things that, unless $n=2$, the canonical ring of $S$ is generated in degree less than or equal to $3$. Theorem~\ref{gen.in.deg2} gives a nice, uniform formula for the number of generators in degree $2$, depending only on the geometric and arithmetic genus of $S$.  Thus a natural question to ask is whether there exists a uniform formula or pattern for the number of generators in degree $3$, similar to the one for generators of degree $2$. Indeed, \cite[Theorem 2.1]{tams} gives us such a formula when $S$ is regular (precisely, the number of extra generators in degree $3$ is $p_g(S)-3$ if $S$ is regular and $n \neq 2$). However, irregular quadruple Galois canonical covers of Types 5.1, 6.1, 7, 8 show that this formula cannot be generalized for an arbitrary $S$, as one can see in Theorem~\ref{canonical.ring} below:

\begin{theorem}\label{canonical.ring}
 Let $\varphi:X \longrightarrow  W$ be a quadruple Galois canonical cover of a surface of minimal degree. Then  the canonical ring of $X$ is generated in degree less than or equal to $3$. Precisely, the canonical ring of $X$ is generated by its part of degree $1$, by $2p_g(X)-4-q(X)$ extra generators in degree $2$, and by $\delta(X)$ extra generators in degree $3$, where
\begin{enumerate}
\item[a)] $\delta(X)=p_g(X)-3$ if $X$ is of Type 1, 2, 3, 4, 5.2, 6.2, 9, 10, 11 or 12. 
\item[b)] $\delta(X)=4$, if $X$ is of Type 5.1 or 6.1;
\item[c)] 
$\delta(X)=5m-1$ if $X$ is of Type 7;
\item[d)] $\delta(X)=9$, if $X$ is of Type 8.1;
\item[e)] 
$\delta(X)=5m$ if $X$ is of Type 8.2.
\end{enumerate}
\end{theorem}

\begin{proof} Let $R$ be the canonical ring of $X$ and $R_n=H^0(nK_X)$ its part of degree $n$. If $X$ is regular (i.e., if $X$ is of Type 1, 2, 3, 4, 9, 10, 11 or 12) the result follows from~\cite[Theorem 2.1]{tams}. The part of the result regarding $R_2$, for any surface $X$, follows from Theorem~\ref{gen.in.deg2}. 

\smallskip Then it only remains to prove that $R$ is generated in degree less than or equal to $3$ and to find the number of extra generators of degree $3$, if $X$ is of Type 5, 6, 7 or 8. 
For this we will study the multiplication maps of global sections of line bundles on $X$
\begin{equation*}
H^0(K_X) \otimes H^0(nK_X)  \overset{\gamma_n}  \longrightarrow H^0((n+1)K_X),
\end{equation*}
when $n \geq 2$. 
The dimension of the cokernel of $\gamma_2$ is the number of linearly independent elements of degree $3$  which are required, together with the elements of $R_1$ and $R_2$, to generate $R_3$. On the other hand, if $\gamma_n$ is surjective for all $n \geq 3$, then $R$ is generated by its elements of degree less than or equal to $3$. 
Now, to study the maps $\gamma_n$, we look at the $\mathcal O_W$--algebra structure of $\varphi_*\mathcal O_X$ and use Lemma~\ref{surject.algstruct} and arguments similar to those used in Section~\ref{Normal.generation.2K_X}.
Recall that by Theorem~\ref{quad12}, $W=\mathbf{F}_{0}$ and that
 the splitting~\eqref{split} of $\varphi_* \mathcal{O}_{X}$ is
\begin{equation}\label{splittings}
\begin{matrix}
\text{if $X$ is of Type 5,}&\mathcal{O}%
_{W}\oplus\mathcal{O}_{W}(-C_{0}-(m+2)f)\oplus\mathcal{O}_{W}(-2C_0)\oplus
\mathcal{O}_{W}(-3C_{0}-(m+2)f);\\
\text{if $X$ is of Type 6,}&\mathcal{O}%
_{W}\oplus\mathcal{O}_{W}(-C_{0}-(m+2)f)\oplus\mathcal{O}_{W}(-2C_0)\oplus
\mathcal{O}_{W}(-3C_{0}-(m+2)f);\\
\text{if $X$ is of Type 7,}&\mathcal{O}%
_{W}\oplus\mathcal{O}_{W}(-(m+1)f)\oplus\mathcal{O}_{W}(-3C_{0}-f)\oplus
\mathcal{O}_{W}(-3C_{0}-(m+2)f)\\
\text{if $X$ is of Type 8,}&\mathcal{O}%
_{W}\oplus\mathcal{O}_{W}(-(m+2)f)\oplus\mathcal{O}_{W}(-3C_{0})\oplus
\mathcal{O}_{W}(-3C_{0}-(m+2)f). 
\end{matrix}
\end{equation}

\noindent We look first at $\gamma_2$. Recall that $H^{0}(\varphi_*K_{X})$ has in all cases one nonzero summand, namely, $H^0(C_0+mf)$. First we consider $X$ of Type 5.1 or 6.1. It follows from~\ref{splittings} that 
the splitting~\eqref{split.pushdown} of $H^{0}(\varphi_*2K_{X})$ has two nonzero summands, namely $H^0(2C_0+2f)$ and $H^0(2f)$, corresponding to $\mathcal O_W$ and $L_2^*$ of~\ref{split}. On the other hand, the splitting~\eqref{split.pushdown} of $H^{0}(\varphi_*3K_{X})$ has four non zero summands, namely, $H^0(3C_0+3f), H^0(2C_{0}), H^0(C_0+3f)$ and 
$H^0(\mathcal O_W)$. Then 
Lemma~\ref{surject.algstruct} implies that $\gamma_2$ is not surjective.  More precisely, since $\gamma_2$ splits in several summands according to the algebra structure of $\varphi_*\mathcal O_{X}$, described in~\cite[Remark 3.1]{GPsmooth},
it follows that the image of $\gamma_2$ is contained in $H^0(3C_0+3f) \oplus H^0(C_0+3f)$. Moreover, $\gamma_2$ surjects onto $H^0(3C_0+3f) \oplus H^0(C_0+3f)$ if the following multiplication maps on $W$
\begin{equation*}
\begin{matrix}
 H^0(C_0+f) \otimes H^0(2C_0+2f) \longrightarrow H^0(3C_0+3f) \\
H^0(C_0+f) \otimes H^0(2f) \longrightarrow H^0(C_0+3f)
\end{matrix}
\end{equation*}
surject. This follows from Lemma~\ref{surjonHirz}, e). Then the dimension of the cokernel of $\gamma_2$ is $h^0(2C_{0})+h^0(\mathcal O_W)=4$.

Arguing similarly we see that in all other Types 5.2, 6.2, 7 and 8, we note that,  from~\ref{splittings}, 
the splitting~\eqref{split.pushdown} of $H^{0}(\varphi_*2K_{X})$ has these nonzero summands:
\begin{enumerate}
 \item[] $H^0(2C_0+2mf) \oplus H^0(C_0+(m-2)f) \oplus H^0(2mf)$, if $X$ is of Type 5.2 or 6.2 (recall that in such a case, $m \geq 2$); 
\item[] $H^0(2C_0+2mf) \oplus H^0(2C_0+(m-1)f)$, if $X$ is of Type 7; 
\item[] $H^0(2C_0+2f)$, if $X$ is of Type 8.1; and 
\item[] $H^0(2C_0+2mf) \oplus H^0(2C_0+(m-2)f)$, if $X$ is of Type 8.2 (recall that in such a case, $m \geq 2$). 
\end{enumerate}
On the other hand, 
since $\gamma_2$ splits in several summands according to the algebra structure of $\varphi_*\mathcal O_{X}$, described in~\cite[Remark 3.1]{GPsmooth},
it follows that the image of $\gamma_2$ is contained in the subspace $B$ of $H^0(\varphi_*3K_{X})$, where $B$ is 
\begin{enumerate}
\item[] $H^0(3C_0+3mf) \oplus H^0(2C_0+(2m-2)f) \oplus H^0(C_0+3mf)$, if $X$ is of Type 5.2 or 6.2 
\item[] $H^0(3C_0+3mf) \oplus H^0(3C_0+(2m-1)f)$, if $X$ is of Type 7; 
\item[] $H^0(3C_0+3f)$, if $X$ is of Type 8.1; and 
\item[] $H^0(3C_0+3mf) \oplus H^0(3C_0+(2m-2)f)$, if $X$ is of Type 8.2 
\end{enumerate}
By Lemma~\ref{surjonHirz}, e), $\gamma_2$ in fact surjects onto $B$ so the cokernel of $\gamma_2$ is 
\begin{enumerate}
\item[] $H^0((2m-2)f)$, if $X$ is of Type 5.2 or 6.2;
\item[] $H^0((3m-1)f) \oplus H^0((2m-2)f)$, if $X$ is of Type 7; 
\item[] $H^0(3C_0) \oplus H^0(3f) \oplus H^0(\mathcal O_W)$, if $X$ is of Type 8.1; and
\item[] $H^0(3mf) \oplus H^0((2m-2)f$, if $X$ is of Type 8.2.
\end{enumerate}
Then the dimension of the cokernel of $\gamma_2$ is 
\begin{enumerate}
\item[]  $2m-1$ if $X$ is of Type 5.2 or 6.2; 
\item[]  $5m-1$  if $X$ is of Type 7;
\item[]  $9$ if $X$ is of Type 8.1; and
\item[] $5m$ if $X$ is of Type 8.2. 
\end{enumerate}
Recall that $2m$ is the degree of $W$, and $W$ is a surface of minimal degree, so $2m=p_g(X)-2$; therefore the number of extra generators of $R$ in degree $3$ is $p_g(X)-3$ if $X$ is of Type 5.2 or 6.2.

\smallskip

Finally
from~\ref{splittings}, 
the splitting~\eqref{split.pushdown} of $H^{0}(\varphi_*nK_{X})$ has four nonzero summands if $n \geq 3$. Then 
the splitting of $\gamma_n$ in summands according to the algebra structure of $\varphi_*\mathcal O_{X}$ and Lemma~\ref{surjonHirz}, e) imply the surjectivity of $\gamma_n$ for all $n \geq 3$. 

\end{proof}

\begin{remark}
{\rm{It is a common phenomenon that, if a graded ring $R$ is generated in degree less than or equal to
$2$, then the Veronese
$2$--subring $R^{\prime}$ of $R$ is generated in degree $1$.  Theorems~\ref{types12},~\ref{types34},~\ref{types5262} and~\ref{canonical.ring} show the existence of several families of surfaces $X$ of general type (both regular and irregular; precisely, surfaces $X$ of Types 1, 2, 3, 4, 5.2 and 6.2) such that the Veronese
$2$--subring of their canonical
ring is generated in degree $1$, despite the fact that their canonical
ring is not generated in degree less than or equal to $2$.}}
\end{remark}


\begin{thebibliography}{Mum90}

\bibitem[Bad01]{Ba}L. B\u adescu, \emph{Algebraic surfaces}, Universitext, Springer--Verlag, 
 New York, 2001. 

\bibitem[BS00]{BS} M. Beltrametti and T. Szemberg, \emph{On higher order embeddings of Calabi--Yau threefolds}, 
Arch. Math. (Basel) {\bf 74} (2000),  221--225. 

\bibitem[Cas06]{Cas} G. Casnati, \emph{Examples of Calabi--Yau threefolds as covers of almost--Fano threefolds}, Geom. Dedicata {\bf 119} (2006), 169--179.

\bibitem[GP98]{GPCalabi} F.J. Gallego and B.P Purnaprajna, \emph{Very ampleness and higher syzygies for Calabi--Yau threefolds}, 
Math. Ann. {\bf 312} (1998), 133--149.

\bibitem[GP01]{rational} F.J. Gallego and B.P Purnaprajna, \emph{Some results on rational surfaces and Fano varieties}, J. Reine Angew. Math. {\bf 538} (2001), 25--55.

\bibitem[GP03]{tams} F.J. Gallego and B.P Purnaprajna, \emph{On the canonical rings of covers of surfaces of minimal degree}, Trans. Amer. Math. Soc. {\bf 355} (2003),  2715--2732

\bibitem[GP07]{GPsing} F.J. Gallego and B.P Purnaprajna, \emph{Classification of quadruple Galois canonical covers II}, J. Algebra {\bf 312} (2007), 798--828.

\bibitem[GP08]{GPsmooth} F.J. Gallego and B.P Purnaprajna, \emph{Classification of quadruple Galois canonical covers I},  Trans. Amer. Math. Soc.  {\bf 360}  (2008),  5489--5507.

\bibitem[Gre82]{Gr} M. Green, \emph{The canonical ring of a variety of general type},
Duke Math. J. {\bf 49} (1982),  1087--1113. 

\bibitem[Hor76]{Ho} E. Horikawa, {\it Algebraic surfaces of general type with
 small $c^2_1$, I}, Ann. of Math. (2) {\bf 104} (1976), 357--387.

 \bibitem[Kon91]{Kon} K. Konno, {\it Algebraic surfaces of general type with $c_1^2
 =3p_g-6$}, Math. Ann. {\bf 290} (1991), 77--107.

\bibitem[Lip69]{Lip} J. Lipman, \emph{Rational singularities, with applications to algebraic surfaces and unique factorization}, Inst. Hautes Études Sci. Publ. Math. {\bf 36} (1969), 195--279.

\bibitem[Mum70]{Mum} D. Mumford, \emph{Varieties defined by quadratic equations}, Corso CIME in Questions on Algebriac Varieties, Rome (1970), 30--100. 

\bibitem[OP95]{OP} K. Oguiso and T. Peternell, \emph{On polarized canonical Calabi--Yau threefolds},  Math. Ann. {\bf 301}  (1995), 237--248. 

\bibitem[Xia90]{Xiao} G. Xiao, \emph{Degree of the bicanonical map of a surface of general type}, 
Amer. J. Math. {\bf 112} (1990), 713--736. 



\end{thebibliography}
\end{document}